\documentclass[a4paper,11pt,leqno]{amsart}

\usepackage{amsmath,amssymb,array,euscript,amsfonts}
\usepackage{graphicx,epsfig,color}
\usepackage[usenames,dvipsnames]{xcolor}
\usepackage{rotating}
\usepackage{bbm}
\usepackage{color}
\usepackage{hyperref}
\def\C{\hbox{\font\dubl=msbm10 scaled 1100 {\dubl C}}}
\def\R{\hbox{\font\dubl=msbm10 scaled 1100 {\dubl R}}}

\def\N{\hbox{\font\dubl=msbm10 scaled 1100 {\dubl N}}}
\def\Z{\hbox{\font\dubl=msbm10 scaled 1100 {\dubl Z}}}

\def\d{\,{\rm{d}}}

\newtheorem{Theorem}{Theorem}

\title[Dirichlet Series with Periodic Coefficients]{Dirichlet Series with Periodic Coefficients and their Value-Distribution Near the Critical Line}

\thanks{{\bf Dedicated to the Memory of Ivan Matveevich Vinogradov at the Occasion of his 130th Birthday}}

\thanks{$^{*}$ The first author is supported by Austrian Science Fund, Project Y-901 and the third author is supported by JSPS KAKENHI Grant Number 18K13400.}
 
\author[Athanasios Sourmelidis, J\"orn Steuding, and Ade Irma Suriajaya]{Athanasios Sourmelidis, J\"orn Steuding, and Ade Irma Suriajaya}

\date{July 2020}

\begin{document}

\begin{abstract}
The class of Dirichlet series associated with a periodic arithmetical function $f$ includes the Riemann zeta-function as well as Dirichlet $L$-functions to residue class characters. We study the value-distribution of these Dirichlet series $L(s;f)$, resp. their analytic continuation in the neighborhood of the critical line (which is the abscissa of symmetry of the related Riemann-type functional equation). In particular, for a fixed complex number $a\neq 0$, we prove for an even or odd periodic $f$ the number of $a$-points of the $\Delta$-factor of the functional equation, prove the existence of the mean-value of the values of $L(s;f)$ taken at these points, show that the ordinates of these $a$-points are uniformly distributed modulo one and apply this to show a discrete universality theorem.
\end{abstract}

\maketitle

{\small \noindent {\sc Keywords:} Dirichlet $L$-functions, Dirichlet series, periodic coefficients, critical line, uniform distribution, universality, Julia line}

\section{Motivation and Statement of the Main Results}

Since the classical works of Euler and Riemann it is known that the distribution of prime numbers is intimately related to the Riemann zeta-function $\zeta(s)$. In general, many topics in analytic number theory are tied to associated Dirichlet series. With his path-breaking articles Ivan Matveevich Vinogradov made outstanding contributions in this direction. For example, he proved the so far widest zero-free region for $\zeta(s)$ and therefore the up to date smallest error term in the prime number theorem \cite{vino2} (independently obtained by Nikolai Michailovitch Korobov); moreover, with his treatment of exponential sums Vinogradov proved that all sufficiently large odd integers satisfy the ternary Goldbach conjecture \cite{vino1} (which was recently extended by Harald Helfgott to all odd integers $\geq 7$). 
\par

In this article we consider generalizations of the Riemann zeta-function and their value-distribution. Given an arithmetical function $f\,:\,\N\to\C$, the associated Dirichlet series is a complex-valued function of a complex variable $s:=\sigma+it$, given by 
$$
L(s;f):=\sum_{n=1}^\infty \dfrac{f(n)}{n^{s}}.
$$
Let $q$ be a positive integer. 
If $f$ is $q$-periodic, that means $f(n+q)=f(n)$ for all positive integers $n$, then the $L(s;f)$ defining Dirichlet series converges absolutely in the right half-plane $\sigma>1$. Moreover, $L(s;f)$ can be analytically continued to the whole complex plane except for at most a simple pole at $s=1$; and $L(s;f)$ is regular at $s=1$ if, and only if the related residue
$$
{\dfrac{1}{q}}\sum_{a\bmod\,q}f(a)
$$
vanishes. 
For technical reasons it is useful to extend $f$ to be defined on $\Z$ by periodicity. In addition, $L(s;f)$ satisfies a functional equation (or {\it identity}), namely,
\begin{eqnarray}\label{feq1}
L(1-s;f)&=&\left(\dfrac{q}{2\pi}\right)^s{\dfrac{\Gamma(s)}{\sqrt{q}}}\left(\mathrm{e}\left(\dfrac{s}{4}\right)L\left(s;f^-\right)+\mathrm{e}\left(-\dfrac{s}{4}\right)L\left(s;f^+\right)\right),
\end{eqnarray}
where $\mathrm{e}(z):=\exp(2\pi iz)$ and $f^{\pm}$ are $q$-periodic arithmetical functions defined by 
\begin{equation}\label{dft}
f^{\pm}(n)={\dfrac{1}{\sqrt{q}}}\sum_{a\bmod\,q}f(a)\mathrm{e}\left(\pm \dfrac{an}{q}\right)
\end{equation}
and may be interpreted as discrete Fourier transforms of $f_{\pm}$ defined by $f_+=f$ and $f_-(n)=f(-n)$. 
All these properties follow from similar properties of the Hurwitz zeta-function $\zeta(s;\alpha)=\sum_{m\geq 0}(m+\alpha)^{-s}$ (with a real parameter $\alpha\in(0,1]$) and a straightforward representation of such $L(s;f)$ as a sum of Hurwitz zeta-functions with rational parameters $\alpha=a/q$. 
This and more were first discovered by Walter Schnee \cite{schnee} ninety years ago. 
\par 

We are concerned with {\it even} and {\it odd} $q$-periodic functions $f$, i.e., $f=\delta f_-$ with $\delta=+1$ for {\it even} $f$ and $\delta=-1$ for {\it odd} $f$. In this case (\ref{feq1}) can be rewritten as
\begin{eqnarray}\label{feq2}
L(s;f)&=&\Delta(s;f)L\left(1-s;f^+\right),
\end{eqnarray}
where
\begin{eqnarray}\label{delta}
\Delta(s;f):=-i\left({\dfrac{q}{2\pi}}\right)^{1-s}{\dfrac{\Gamma(1-s)}{\sqrt{q}}}\left(\mathrm{e}\left(\dfrac{s}{4}\right)-\delta \mathrm{e}\left(-\dfrac{s}{4}\right)\right).
\end{eqnarray}
This setting includes the case of Dirichlet $L$-functions associated with (not necessarily primitive) residue class characters $\chi\bmod\,q$ (see \cite{apostol} for details). Some of our results generalize previous ones for the Riemann zeta-function \cite{Kalp,korolev,steusur}; however, the findings concerning uniform distribution and universality are new (although the latter property has been considered in this context earlier by Maxim Korolev and Antanas Laurin\v cikas \cite{kola} for a special case). Our approach is inspired by rather general though deep theorems from the theory of functions due to \'Emile Picard, Gaston Julia and Rolf Nevanlinna, and our analysis relies mainly on the functional equation \eqref{feq2}.
\par

In 1879, Picard \cite{picard} proved that if an analytic function $f$ has an essential singularity at a point $\omega$, then $f(s)$ takes all possible complex values with at most a single exception (infinitely often) in every neighbourhood of $\omega$. And in 1919, a little more than one hundred years ago, Julia \cite{julia} refined Picard's great theorem by showing that one can even add a restriction on the angle at $\omega$ to lie in an arbitrarily small cone (see also \cite[\S 12.4]{burckel}). 
For Dirichlet series appearing in number theory, however, it is more natural to consider so-called {\it Julia lines} rather than Julia directions. 
For this and further references we refer to \cite{christ,steudi,steusur}.
\par 

Given a complex number $a$, the solutions of the equation
$$
\Delta\left(s;f\right)=a
$$
are called the {\it $a$-points} of $\Delta\left(s;f\right)$, and we denote them as $\delta_a:=\beta_a+i\gamma_a$. 
We shall show that for any fixed $a\neq 0$ {\it most} of the $a$-points are clustered around the critical line $1/2+i\R$ which, therefore, is a Julia line for $\Delta\left(s;f\right)$. Moreover, we prove that the mean of the values $L(\delta_a;f)$ exists and that every complex number appears as such a limit. This indicates an interesting link between the distribution of $a$-points of $\Delta(s;f)$ and the values taken by $L(s;f)$. For the case of the Riemann zeta-function and values $a$ from the unit circle these $a$-points have been studied by Justas Kalpokas and the second author \cite{Kalp} as well as by Kalpokas, Maxim Korolev and the second author \cite{korolev}, respectively. In this situation the $\exp(2i\phi)$-points correspond to intersections of the curve $t\mapsto \zeta\left({\frac{1}{ 2}}+it\right)$ with straight lines $\exp(i\phi)\R$ through the origin. 

\begin{Theorem}\label{julia}
Let $f$ be an even or odd $q$-periodic function. Then, the function $s\mapsto \Delta\left(s;f\right)$ is a meromorphic function with two exceptional values $0$ and $\infty$ in the sense of Nevanlinna theory and the critical line $1/2+i\R$ is a Julia line for $\Delta(s;f)$. Moreover, given a complex number $a\neq 0$, the number $N_a(T;f)$ of $a$-points $\delta_a= \beta_a+ i\gamma_a$ of $\Delta\left(s;f\right)$ satisfying $0<\beta_a<1, 0<\gamma_a<T$ is asymptotically equal to
\begin{equation}\label{na}
N_a(T;f)=\dfrac{T}{2\pi}\log{\dfrac{qT}{2\pi e}}+O(\log T).
\end{equation}
\end{Theorem}

\noindent The error term here and all error terms in the sequel depend on $a$ and $f$. The condition $0<\beta_a<1$ in the above theorem can be relaxed to any open interval centered at $1/2$. This can be easily explained using an asymptotic equation for $\Delta(s;f)$ (appearing in the next section).

\begin{Theorem}\label{mean}
Let $f$ be an even or odd $q$-periodic function and let $\delta_a=\beta_a+i\gamma_a$ denote the $a$-points of $\Delta(s;f)$. Then, for every complex number $a\neq 0$ and any $\epsilon>0$, 
\begin{eqnarray*}
\mathop{\sum_{0<\gamma_a<T}}_{ 0<\beta_a<1}L(\delta_a;f)&=&\left(f(1)+af^+(1)\right)\dfrac{T}{2\pi}\log{\dfrac{qT}{2\pi e}}+O_\epsilon\left(T^{1/2+\epsilon}\right).
\end{eqnarray*}
\end{Theorem}

The horizontal distribution of the $a$-points of $\Delta(s;f)$ is quite regular as we will see in the next section. We want to provide some additional information regarding their vertical distribution, since such questions arise quite often in zeta-function theory. For example, in the case of $\zeta(s)$, Edmund Landau \cite{land} proved the explicit formula
\begin{equation}\label{Lan}
\sum\limits_{0<\gamma\leq T}x^\rho=-\dfrac{T}{2\pi}\Lambda(x)+O_x(\log T),
\end{equation}
valid for any $T>1$ and any positive number $x>1$, where the $\rho:=\beta+i\gamma$ denote the non-trivial zeros of $\zeta(s)$ and $\Lambda(x)$ is the von Mangoldt function extended to the whole real line by setting it to be zero if $x$ is not a positive integer. The second author \cite{steu1} proved a similar result in the case of $a$-points of $\zeta(s)$, namely, for $a\neq1$,
\begin{equation}\label{St}\sum\limits_{0<\gamma_a\leq T}x^{\rho_a}=\dfrac{T}{2\pi}\left(a(x)-x\Lambda\left(\frac{1}{x}\right)\right)+O_{x,\epsilon}\left(T^{1/2+\epsilon}\right),
\end{equation}
where $\rho_a:=\beta_a+i\gamma_a$ denotes the so called {\it non-trivial} $a$-points of $\zeta(s)$ and $a(x)$ is a certain computable arithmetical function extended to the whole real line. The authors could not succeed in proving a comparable asymptotic formula in the case of the $\delta_a$ but an upper bound for the corresponding sum with respect to the $a$-points of $\Delta$ which is remarkably small compared with \eqref{Lan} or \eqref{St}. 

In order to formulate this result we recall some standard notation. The number $\lfloor x\rfloor$ denotes the greatest integer which is smaller than or equal to the real number $x$ and $\mathbbm{1}_S$ denotes the characteristic function of a set $S$.

\begin{Theorem}\label{uniform}
Let $a\neq0$ be a fixed complex number. 
Then there exists a constant $c>0$, depending on $a$, such that, for any $0<x\neq1$ and any $T,T'$ satisfying $\max\left\{1,4\pi/q\right\}\leq T<T+1\leq T'\leq2T$, we have that
\begin{eqnarray*}
\sum_{\substack{T<\gamma_a<T'}} x^{\delta_a} &\ll& x^{1/2}\left(x^{c/\log \frac{qT}{2\pi}}+x^{-c/\log \frac{qT}{2\pi}}\right)\left(|\log x|+\log T+\dfrac{\log T}{|\log x|}\right) \\
&&+~ \mathbbm{1}_{(1,+\infty)}(x) E_1(x,a,T)-\mathbbm{1}_{(0,1)}(x)E_2(x,a,T),
\end{eqnarray*}
where $E_1(x,a,T)$ and $E_2(x,a,T)$ are always zero, unless there is a positive integer $$j\leq\left\lfloor\dfrac{2\log\frac{qT}{4\pi}}{\log30}\right\rfloor$$
such that $x^{1/j}\in\left[qT/(4\pi),5qT/(4\pi)\right]$ or $x^{-1/j}\in\left[qT/(4\pi),5qT/(4\pi)\right]$, in which case we have that
\begin{eqnarray*}
E_1(x,a,T)&:=&\dfrac{x^{1/2}\log x}{j^{3/2}}\left(\dfrac{x^{1/(2j)}}{|a|^{j}}+\dfrac{x^{c/\log \frac{qT}{2\pi}}T^{1/2}}{30^{j}}\right),\\
E_2(x,a,T)&:=&\dfrac{x^{1/2}\log x}{j^{3/2}}\left({x^{-1/(2j)}}|a|^j+\dfrac{x^{-c/\log \frac{qT}{2\pi}}T^{1/2}}{30^{j}}\right),
\end{eqnarray*}
respectively. 
In particular, if $x\neq1$ is such that $4\pi/(qT)\leq x\leq qT/(4\pi)$, then
$$ \sum_{\substack{T<\gamma_a<T'}} x^{\delta_a}\ll x^{1/2}\log T\left(1+\dfrac{1}{|\log x|}\right).
$$
\end{Theorem}

A sequence of real numbers $(x_n)_{n\in\mathbb{N}}$ is called {\it uniformly distributed modulo 1} if
$$
\lim\limits_{N\to\infty}\dfrac{1}{N}\sharp\left\{1\leq n\leq N:\lbrace x_n\rbrace\in[a,b]\right\}=b-a
$$
for any real numbers $0\leq a\leq b\leq1$, where $\lbrace x\rbrace:=x-\lfloor x\rfloor$ for $x\geq0$. 
Hans Rademacher \cite{rader} employed Landau's formula \eqref{Lan} to prove, under the Riemann Hypothesis, that the sequence $(\alpha\gamma)_{\gamma>0}$ is uniformly distributed modulo $1$, where $\alpha\neq0$ is a real number and $(\gamma)_{\gamma>0}$ is the sequence of ordinates of the non-trivial zeros of $\zeta(s)$ in ascending order (and counted with multiplicities). 
Peter Elliott \cite{elliot} and (independently) Edmund Hlawka \cite{hlawka} proved the uniform distribution of $(\alpha\gamma)_{\gamma>0}$ unconditionally. 
Respectively, the second author \cite{steu1} proved the uniform distribution of the sequence $(\alpha\gamma_a)_{\gamma_a>0}$ , where $(\gamma_a)_{\gamma_a>0}$ is the sequence of ordinates of the non-trivial $a$-points of $\zeta(s)$ in ascending order (counting multiplicities). 
The common feature of both statements is that they follow from \eqref{Lan} and \eqref{St} and known results on the {\it clustering} of the non-trivial zeros (or $a$-points) of $\zeta(s)$ around the critical line. This finally leads to
\begin{equation}\label{ig}
\sum\limits_{0<\gamma<T}x^{i\gamma}=o_x(T\log T)=\sum\limits_{0<\gamma_a<T}x^{i\gamma_a}
\end{equation}
and uniform distribution follows from a well-known criterion of Hermann Weyl. It is noteworthy that, for fixed $x$, the best possible bounds known for the sums in \eqref{ig} are of order $T$. 
Compared to these results, the bounds given in the following theorem for the case of $a$-points of $\Delta(s;f)$ are very different. This is probably not unexpected since $\Delta(s;f)$ is not a zeta-function to begin with and the clustering of its $a$-points around the critical line has a simple explanation as we shall see in Sections \ref{juliaa} and \ref{uniform}.

\begin{Theorem}\label{u.d.}
Let $0<\gamma_a^{(1)}\leq\gamma_a^{{(2)}}\leq\gamma_a^{{(3)}}\leq\dots$ be the sequence of ordinates of all $a$-points of $\Delta(s;f)$ (counted with multiplicities) and having real part in $(0,1)$, where $f$ is an odd or even $q$-periodic function. Then, for any integer $N\geq0$ and any number $x\neq1$ satisfying $4\pi/(qN)\leq x\leq qN/(4\pi)$, we have 
\begin{equation*}
\sum\limits_{N< n\leq 2N}x^{i\gamma_a^{{(n)}}}\ll\left(\dfrac{1}{|\log x|}+|\log x|\right)\log N.
\end{equation*}
In particular, the sequence of numbers $\alpha\gamma_a^{{(n)}}$, $n\in\mathbb{N}$, is uniformly distributed modulo $1$ for any real number $\alpha\neq0$.
\end{Theorem}

Robert Spira \cite{spira} showed that the $\Delta$-function appearing in the functional equation for the Riemann zeta-function $\zeta(s)=L(s;\mathbf{1})$ (with $\mathbf{1}$ being the function constant $1$) satisfies $\vert \Delta\left(s;\mathbf{1}\right)\vert <1$ for $1/2<\sigma<1$ and $t\geq10$. This implies along with the identity $\Delta\left(s;\mathbf{1}\right)\Delta\left(1-s;\mathbf{1}\right)=1$ that, if $a$ is a complex number from the unit circle, then all but finitely many $a$-points of $\Delta\left(s;\mathbf{1}\right)$ have real part equal to $1/2$. Therefore, our theorems generalize some results from \cite{Kalp},
as well as recent results by Korolev \& Laurin\v cikas \cite{kola} on the uniform distribution of Gram points (the $1$-points of $\Delta\left(s;\mathbf{1}\right)$). This last article motivates us to prove the following {\it joint universality theorem}:

\begin{Theorem}\label{universal}
Let $f$ be an even or odd $q$-periodic function and denote the ordinates of the $a$-points of $\Delta(s;f)$ in ascending order by $\gamma_a^{(n)}$. 
Let also $\chi_1,\chi_2,\dots,\chi_J$ be non-equivalent Dirichlet characters and $\psi\neq\mathbf{0}$ be an $r$-periodic arithmetical function with $r\neq2$. 
Then, for any compact set with connected complement $K$ inside the strip $1/2<\sigma<1$, any $g_1,g_2,\dots,g_J,h$ continuous non-vanishing functions on $K$ which are analytic in its interior, any real numbers $z>0$ and $\xi_p$, indexed by the primes $p\leq z$, and any $\varepsilon>0$,
\begin{equation}\label{joint}
\liminf\limits_{N\to\infty}\dfrac{1}{N}\sharp\left\{1\leq n\leq N:\begin{array}{c}\max\limits_{1\leq j\leq J}\max\limits_{s\in K}\left|L\left(s+i\gamma_a^{{(n)}};\chi_j\right)-g_j(s)\right|<\varepsilon,\\ \& \ \max\limits_{p\leq z}\left\|\gamma_a^{{(n)}}\dfrac{\log p}{2\pi}-\xi_p\right\|<\varepsilon
\end{array}\right\}>0
\end{equation}
and
\begin{equation}\label{periodic}
\liminf\limits_{N\to\infty}\dfrac{1}{N}{\sharp\left\{1\leq n\leq N:\max\limits_{s\in K}\left|L\left(s+i\gamma_a^{{(n)}};\psi\right)-h(s)\right|<\varepsilon\right\}}>0,
\end{equation}
where $\sharp A$ denotes the cardinality of a set $A\subseteq\mathbb{R}$ and $\|x\|:=\min_{m\in\mathbb{Z}}|m-x|$.
In addition, if $\psi$ has period $r\geq3$ and is not a multiple of a Dirichlet character mod $r$, then $h$ is allowed to have zeros in $K$.
In the case when $r=2$, \eqref{periodic} holds for any compact set with connected complement $K$ inside the open set
$$D_0:=\left\{s\in\mathbb{C}:\dfrac{1}{2}<\sigma<1\right\}\setminus\left\{\log\left(1-\dfrac{\psi(2)}{\psi(1)}\right)+2k\pi i:k\in\mathbb{Z}\right\},$$
where $\log$ is the principal logarithm.
\end{Theorem}

\noindent This is a discrete version of Sergei Voronin's classical universality theorem for the Riemann zeta-function \cite{voronin}: 
\begin{eqnarray}\label{univer}
\liminf\limits_{T\to\infty}\dfrac{1}{T}{\sharp\left\{\tau\in[0,T]:\max\limits_{s\in K}\left|\zeta(s+i\tau)-h(s)\right|<\varepsilon\right\}}>0
\end{eqnarray}
(resp. its simultaneous version \cite{voro2} for a family of Dirichlet $L$-functions which is often called {\it joint universality}). As a matter of fact, Voronin proved that there exists some real $\tau$ such that $\zeta(s+i\tau)$ is $\varepsilon$-close to $h(s)$ when $s$ ranges in a disc of center $3/4+it_0$ and radius $r<1/4$ (see for example \cite[Chapter VII]{karat}). A few years later Axel Reich \cite{reich} and (independently) Bhaskar Bagchi \cite{bagchi} obtained \eqref{univer} (which is implicit in Voronin's work) and provided also a discrete version with respect to arithmetic progressions
\begin{eqnarray*}\label{duniver}
\liminf\limits_{T\to\infty}\dfrac{1}{N}{\sharp\left\{1\leq n\leq N:\max\limits_{s\in K}\left|\zeta\left(s+id n\right)-h(s)\right|<\varepsilon\right\}}>0,
\end{eqnarray*}
where $d$ is a fixed non-zero real number. Theorem \ref{universal} is of similar nature, where the shifts are ordinates of $a$-points of $\Delta(s;f)$.

\section{Proof of Theorem \ref{julia}}\label{juliaa}

We first remark that $\Delta(s;f)$ depends only on the period $q$ of $f$ and $\delta=\pm1$ which determines whether $f$ is an even or odd function. We also observe that if $f$ is an even or odd $q$-periodic arithmetical function, then so are $f^\pm$ and $\overline{f}$ (the complex conjugate of $f$). In the rest of this paper, we will use the simplified notation 
\begin{eqnarray}\label{delta1}
\Delta(s):=\Delta(s;f)=\Delta\left(s;\overline{f}\right)=\Delta\left(s;f^+\right).
\end{eqnarray}
This function $\Delta(s)$ is the product of an exponential function, the Gamma-function, and a trigonometric function. 
It is well-known that $\Gamma(z)^{-1}$ is an entire function with only simple zeros at $z=-n$ for $n\in\N_0:=\N\cup\{0\}$. 
Hence, by (\ref{delta}), $\Delta(s)$ is regular except for simple poles at the positive odd integers (if $\delta=+1$), respectively for the positive even integers (if $\delta=-1$); moreover, $\Delta(s)$ vanishes exactly for the non-positive even integers (if $\delta=+1$), respectively for the negative odd integers. One can show by an application of Rouch\'e's theorem (as Levinson did in \cite{levi} for $\Delta(s;\mathbf{1})$, the case of the Riemann zeta-function) that there are a few $a$-points of $\Delta(s)$ in a neighbourhood of the real line; their count inside a strip $\left\{x+iy:|x|\leq r\text{ and }|y|\leq1\right\}$ is $O(r)$ as $r\to\infty$. In comparison with Formula \eqref{na} it follows (along the lines of \cite[Section 7.3]{steudi}) that $0$ and $\infty$ are exceptional values of $\Delta(s)$ in the sense of Nevanlinna theory. It remains to prove the Riemann-von Mangoldt-type formula \eqref{na}. 

It is an easy consequence of Stirling's formula,
\begin{equation} \label{stirling}
\Gamma(\sigma+it)=\sqrt{2\pi}t^{\sigma+it-1/2}\exp\left(-{\dfrac{\pi t}{ 2}}-it+{\dfrac{\pi i}{ 2}}\left(\sigma-\dfrac{1}{2}\right)\right)\left(1+O\left(\dfrac{1}{t}\right)\right),
\end{equation}
and the reflection principle
\begin{equation}\label{reflection}\Gamma(\sigma-it)=\overline{\Gamma(\sigma+it)}
\end{equation} 
both valid uniformly for $t\geq 1$ and $\sigma$ from any strip of bounded width, and (\ref{delta}) that
\begin{equation}\label{pow}
\Delta(\sigma+it)=\delta\left({\dfrac{qt}{ 2\pi}}\right)^{1/2-(\sigma+it)}\exp\left(i\left(t+\dfrac{\pi}{4}\right)\right)\left(1+O\left(\dfrac{1}{t}\right)\right)
\end{equation}
as $t\to+\infty$. 
Hence, $\Delta(\sigma+it)$ tends for $\sigma<1/2$ to infinity and for $\sigma>1/2$ to zero as $t\to +\infty$. 
Thus, the critical line $1/2+i\R$ divides the upper half-plane into two domains where the limit $\lim_{t\to+\infty}\Delta(\sigma+it)$ exists in the compactified plane $\C\cup\{\infty\}$; on the critical line however the limit does not exist. The behaviour in the lower half-plane is ruled by conjugation, 
\begin{eqnarray}\label{refle}
\Delta(\sigma-it) = \delta\,\overline{\Delta(\sigma+it)} &\text{ and }&
\Delta'(\sigma-it) = \delta\,\overline{\Delta'(\sigma+it)},
\end{eqnarray}
as follows from \eqref{delta}, \eqref{reflection} and Cauchy's integral formula. Near the boundary (critical) line $1/2+i\R$, however, the distribution of values is rather different. As a matter of fact, $\Delta(s)$ takes every complex value $a\neq 0$ infinitely often there. Writing $a=\Delta(\delta_a)=\vert a\vert \exp(i\phi)$ with an $a$-point $\delta_a=\beta_a+i\gamma_a$ of $\Delta$ and comparing with (\ref{pow}) implies that 
\begin{eqnarray}\label{amod}
\vert a\vert &=& \left({\dfrac{q\gamma_a}{ 2\pi}}\right)^{1/2-\beta_a} \left(1+O\left(\dfrac{1}{\gamma_a}\right)\right),\\
\phi &\equiv& \gamma_a\log {\dfrac{2\pi e}{ q\gamma_a}}+{\dfrac{\pi}{ 4}}+\dfrac{(1-\delta)\pi}{2}+O\left(\dfrac{1}{\gamma_a}\right)\ \bmod\, 2\pi.\nonumber
\end{eqnarray}
This shows that $\beta_a\to 1/2$ as $\gamma_a\to+\infty$ (and explains a remark from the introduction). In particular, there exists a real number $t_a>0$, depending only on $a,\delta$ and $q$, such that $t_a$ is not an ordinate of any $a$-point of $\Delta$ and $\beta_a\in(0,1)$ if, and only if, $\gamma_a\geq t_a$; we can actually choose $t_a$ here such that $\beta_a$ is included in any open interval centered at $1/2$, but the way we define it here yields, for instance, $N_a(T;f)=\sharp\lbrace \gamma_a:t_a<\gamma_a<T\rbrace$.

Before showing \eqref{na}, we use an argument similar to \cite[\S 9.2]{tit} to prove
\begin{align}\label{app1}
N_a(T+1;f)-N_a(T;f)\ll\log T
\end{align}
for any $T\geq t_a+3$. Indeed, if $n(r)$ denotes the number of $a$-points of $\Delta(s)$ in the disc with center $2+iT$ and radius $r$, then
$$
N_a(T+1;f)-N_a(T;f)\leq n\left(\sqrt{5}\right)\ll\int_0^3\dfrac{n(r)}{r}\d r.
$$
It follows from Jensen's formula (see for example \cite[\S 3.61]{tit1}) and \eqref{pow} that
\begin{eqnarray*}
\int_0^3\dfrac{n(r)}{r}\d r&=&\dfrac{1}{2\pi}\int_0^{2\pi}\log\left|\Delta\left(2+iT+3e^{i\theta}\right)-a\right|\d\theta-\log\left|\Delta\left(2+iT\right)-a\right|\\&\ll&\log T
\end{eqnarray*}
and thus \eqref{app1} holds.

To prove \eqref{na} we apply the argument principle to the function $\Delta(s)-a$ and integrate over the counterclockwise oriented rectangle ${\mathcal C}$ with vertices $-1+it_a$, $2+it_a$, $2+iT$ and $-1+iT$. This gives
\begin{equation}\label{arg_prin}
2\pi i N_a\left(T;f\right) =\int_{\mathcal{C}} {\dfrac{\Delta'(s)}{ \Delta(s)-a}} \d s;
\end{equation}
since all poles of $\Delta(s)$ lie on the real line, they do not affect here. In addition, all zeros lie outside ${\mathcal C}$; hence, we may rewrite the integrand as
\begin{equation}\label{bibo}
{\dfrac{\Delta'(s)}{\Delta(s)-a}}={\dfrac{\Delta'}{ \Delta}}(s)\cdot \frac{1}{1 - a/\Delta(s)}
\end{equation}
or
\begin{equation}\label{bibo2}
{\dfrac{\Delta'(s)}{\Delta(s)-a}}={\dfrac{\Delta'(s)}{-a}}\cdot \frac{1}{1 - \Delta(s)/a}.
\end{equation}
Now, taking into account (\ref{pow}) in combination with another form of Stirling's formula,
\begin{equation} \label{log_der}
{\dfrac{\Delta' }{ \Delta}}(\sigma+it) = -\log{\dfrac{qt}{ 2\pi}} + O\left(\dfrac{1}{t}\right),
\end{equation}
which is also valid for $t\geq t_a>0$ and $\sigma$ from any strip of bounded width, we obtain, for $\epsilon>0$ and $\sigma\geq1/2+\epsilon$,
\begin{equation}\label{one}
{\dfrac{\Delta'(s)}{\Delta(s)-a}}\ll_\epsilon t^{1/2-\sigma}\log t.
\end{equation}
Similarly, we have, for $\epsilon>0$ and $\sigma\leq1/2-\epsilon$, that
\begin{eqnarray}\label{two}
{\dfrac{\Delta'(s) }{ \Delta(s)-a}}&=&{\dfrac{\Delta'}{\Delta}}(s)\cdot \left(1+\sum_{j\geq 1}\left({\dfrac{a}{ \Delta(s)}}\right)^j\right)\\
&=& -\log{\dfrac{qt}{ 2\pi}} + O_\epsilon\left(t^{-1}+t^{\sigma-1/2}\log t\right)\nonumber;
\end{eqnarray}
here (\ref{pow}), obviously, allows us to expand the second factor into a geometric series. 

Next, it follows from \eqref{one} that the contribution of the integral over the right vertical segment is negligible: 
$$
\int_{2+it_a}^{2+iT}{\dfrac{\Delta'(s)}{ \Delta(s)-a}} \d s \ll 1.
$$
Using expression \eqref{two}, however, the contribution of the integral over the left vertical segment leads to
\begin{eqnarray*}
\int_{-1+iT}^{-1+it_a}{\dfrac{\Delta'(s)}{\Delta(s)-a}} \d s
&=& -i\int_{t_a}^T\left(-\log{\dfrac{qt}{ 2\pi}}+O\left(\dfrac{1}{t}\right)\right)\d t+O(\log T)\\
&=&i T\log{\dfrac{qT}{ 2\pi e}}+O(\log T).
\end{eqnarray*}
It remains to estimate the integrals on the horizontal segments of $\mathcal{C}$. The lower one is trivially bounded
$$
\int_{-1+it_a}^{2+it_a}{\dfrac{\Delta'(s)}{\Delta(s)-a}} \d s\ll1,
$$
while for the upper one may use the following truncated partial fraction decomposition
\begin{align}\label{aprox}
{\dfrac{\Delta'(s)}{ \Delta(s)-a}}=\sum_{\vert t-\gamma_a\vert\leq 1}\frac1{s-\delta_a}+O(\log t),
\end{align}
which is valid for $\sigma\in[-1,2]$ and $t\geq t_a>0$. We will show this after finishing the proof of Theorem \ref{julia}. 

In view of (\ref{aprox})
$$
\int_{2+iT}^{-1+iT}{\dfrac{\Delta'(s)}{\Delta(s)-a}} \d s=\sum_{\vert T-\gamma_a\vert\leq 1}\int_{2+iT}^{-1+iT}\frac{\d s}{s-\delta_a}+O(\log T).
$$
By the calculus of the residues we obtain that
$$
\int_{2+iT}^{-1+iT}\frac{\d s}{s-\delta_a}=\left\{\int_{2+iT}^{2+i(T+2)}+\int_{2+i(T+2)}^{-1+i(T+2)}+\int_{-1+i(T+2)}^{-1+iT}\right\}\frac{\d s}{s-\delta_a}-2\pi i R(\delta_a),
$$
where $R(\delta_a)$ is $1$ or $0$ depending on whether $\delta_a$ lies inside the rectangle described above or not. Recall that $\beta_a\in(0,1)$. Hence, for every $a$-point with $\left|T-\gamma_a\right|\leq1$,
\begin{eqnarray*}
\int_{2+iT}^{-1+iT}\frac{\d s}{s-\delta_a}&\ll&\int_T^{T+2}\dfrac{\d t}{\left|2-\beta_a+i\left(t-\gamma_a\right)\right|}+\int_{-1}^{2}\dfrac{\d \sigma}{\left|\sigma-\beta_a+i\left(T+2-\gamma_a\right)\right|}\\
&&+ \int_T^{T+2}\dfrac{\d t}{\left|-1-\beta_a+i\left(t-\gamma_a\right)\right|}+1\\
&\ll&1.
\end{eqnarray*}
In combination with inequality \eqref{app1}, we obtain 
\begin{align*}
\int_{2+iT}^{-1+iT}{\dfrac{\Delta'(s)}{\Delta(s)-a}} \d s\ll\sum_{\vert T-\gamma_a\vert\leq 1}1+O(\log T)\ll\log T.
\end{align*}
Finally, we arrive at
$$
\int_{\mathcal{C}} {\dfrac{\Delta'(s)}{\Delta(s)-a}} \d s=i T\log{\dfrac{qT}{ 2\pi e}}+O(\log T).
$$
Substituting this into \eqref{arg_prin} finishes the proof of \eqref{na}.

It remains to show \eqref{aprox}. For this purpose, we apply Jacques Hadamard's theory of functions of finite order (see \cite[\S 5.3]{palka}); our reasoning is similar to the case of the Riemann zeta-function (see \cite[\S 9.6]{tit}).

As already mentioned, $\Delta(s)$ is analytic except for simple poles at some positive integers. 
Thus, 
$$
F(s):=(\Delta(s)-a)\cdot \Gamma(1-s)^{-1} 
$$
defines an entire function. By Stirling's formula it follows that $F$ is entire and of order one. Hence, Hadamard's factorization theorem implies the product representation
$$
F(s)=\exp(A+Bs)\prod_{\delta_a}\left(1-\dfrac{s}{ \delta_a}\right)\exp\left(\dfrac{s}{ \delta_a}\right),
$$
where $A$ and $B$ are certain complex constants and the product is taken over {\it all} zeros $\delta_a$ of $F(s)$. Taking the logarithmic derivative, we deduce
$$
\dfrac{F'}{ F}(s)=B+\sum_{\delta_a}\left(\frac1{s-\delta_a}+\frac1{\delta_a}\right).
$$
Since 
$$
\dfrac{F'}{ F}(s)=\dfrac{\Delta'(s)}{ \Delta(s)-a}+\dfrac{\Gamma'} {\Gamma}(1-s)
$$
and $\dfrac{\Gamma'}{ \Gamma}(1-s)\ll \log t$ (also a consequence of Stirling's formula), we have
$$
\dfrac{\Delta'(s)}{\Delta(s)-a}=\sum_{\delta_a}\left(\frac1{s-\delta_a}+\frac1{\delta_a}\right)+O(\log t).
$$
Setting $s=2+it$ and subtracting from the latter expression
$$
\dfrac{\Delta'(2+it)}{ \Delta(2+it)-a},
$$
which is $O(1)$ by \eqref{one}, we obtain
$$
\dfrac{\Delta'(s)}{ \Delta(s)-a}=\sum_{\delta_a}\left(\frac1{s-\delta_a}-\dfrac{1}{ 2+it-\delta_a}\right)+O(\log t).
$$
In view of \eqref{app1}, it follows that
$$
\sum_{\vert t-\gamma_a\vert\leq 1}\frac1{2+it-\delta_a}\ll \sum_{\vert t-\gamma_a\vert\leq 1}1=N_a(t+1;f)-N_a(t-1;f)\ll \log t.
$$
A short computation shows that, for any positive integer $n$,
$$
\sum_{t+n<\gamma_a\leq t+n+1}\left(\frac1{s-\delta_a}-\dfrac{1}{ 2+it-\delta_a}\right)\ll \sum_{t+n<\gamma_a\leq t+n+1}\dfrac {1}{ n^2}\ll\dfrac {\log (t+n)}{ n^2}.
$$
And since $\sum_{n\geq 1}\log(t+n)/n^2\ll \log t$, it follows that
$$
\sum_{\gamma_a>t+1}\left(\frac1{s-\delta_a}-\dfrac{1}{ 2+it-\delta_a}\right)\ll \log t;
$$
obviously, we can bound the sum over the $a$-points $\delta_a$ satisfying $\gamma_a<t-1$ similarly by the same bound. Consequently, the contribution of the $a$-points distant from $s$ are negligible. This implies (\ref{aprox}) and concludes the proof of Theorem \ref{julia}.

\section{Proof of Theorem \ref{mean}}

In view of \eqref{na}, it follows that for any given $T_0>0$, there exists $T\in[T_0,T_0+1)$ such that
\begin{equation}\label{condi}
\min_{\delta_a}\vert T-\gamma_a\vert\gg\dfrac{1}{\log T},
\end{equation}
where the minimum is taken over all $a$-points $\delta_a$. In the proof of Theorem \ref{julia} we have observed that there exists a real number $t_a>0$ such that $\beta_a\in(0,1)$ if and only if $\gamma_a>t_a$. Thus integrating over the counterclockwise oriented rectangle ${\mathcal C}$ with vertices $2+it_a, 2+iT, -\epsilon+iT, -\epsilon+it_a$, where $\epsilon>0$, we obtain
\begin{eqnarray}\label{i1234}
\lefteqn{\mathop{\sum_{0<\gamma_a<T}}_{ 0<\beta_a<1}L(\delta_a;f)}\nonumber\\
&=&\frac{1}{2\pi i}\left\{\int_{2+it_a}^{2+iT}\hspace*{-1pt}+\hspace*{-1pt}\int_{2+iT}^{-\epsilon+iT}\hspace*{-1pt}+\hspace*{-1pt}\int_{-\epsilon+iT}^{-\epsilon+it_a}\hspace*{-1pt}+\hspace*{-1pt}\int_{-\epsilon+it_a}^{2+it_a}\right\}\hspace*{-1pt} \frac{\Delta'(s)}{\Delta(s)-a}L(s;f)\d s\nonumber \\
&=:&\sum_{1\leq j\leq 4}\mathcal{I}_j.
\end{eqnarray}

Since the Dirichlet series coefficients $f(n)$ are bounded, we have $L(s;f)\ll_\epsilon 1$ for $\sigma\geq 1+\epsilon$. In view of the functional equation \eqref{feq2} and the asymptotic formula (\ref{pow}) it follows from the Phragm\'en-Lindel\"of principle (which is a kind of maximum principle for unbounded domains) that
\begin{equation}\label{ioi}
\begin{array}{cccc}
L(\sigma+it;f)&\ll_\epsilon& 1+t^{(1-\sigma)/2+\epsilon}&\quad\mbox{for}\quad \sigma\in[0,2],\ t\geq x_0>0,
\end{array}
\end{equation}
and 
\begin{equation}\label{oio}
\begin{array}{cccc}
L(\sigma+it;f)&\ll_\epsilon &t^{1/2-\sigma+\epsilon}&\quad\mbox{for}\quad \sigma\in [-1,0],\ t\geq x_0>0
\end{array}
\end{equation}
(see \cite[\S 5.1]{tit} and \cite[\S 9.41]{tit1} for the case of $\zeta(s)$; the generalization to $L(s;f)$ is straightforward).

We begin with the vertical integrals in (\ref{i1234}). Since the range of integration of $\mathcal{I}_1$ lies in the half-plane $\sigma>1/2$, it follows from (\ref{one}) in combination with \eqref{ioi} that
\begin{eqnarray}\label{i1}
{\mathcal I}_1 \ll \int_{t_a}^{T} t^{-3/2}\log t \d t\ll 1.
\end{eqnarray}
The integral $\mathcal{I}_3$ is over a line segment in $\sigma<1/2$. It follows from the functional equation \eqref{feq2} and relation \eqref{delta1} that
\begin{equation}\label{neu}
L(1-s;f)=\Delta(1-s)L\left(s;f^+\right)=\Delta(1-s)\Delta(s)L\left(1-s;\left(f^{+}\right)^+\right).
\end{equation}
Here we observe an instance of the Fourier inversion formula, namely $\left(f^+\right)^+=\delta f$. In order to see that we compute via (\ref{dft}) that 
\begin{eqnarray*}
\left(f^{\pm}\right)^+(n)&=&{\dfrac{1}{\sqrt{q}}}\sum_{a\bmod\,q}f^{\pm}(a)\mathrm{e}\left(\dfrac{an}{q}\right)\\
&=&{\dfrac{1}{ q}}\sum_{a,b\bmod\,q}f(b)\mathrm{e}\left(\dfrac{\pm ba}{q}\right)\mathrm{e}\left(\dfrac{an}{q}\right)\\
&=&{\dfrac{1}{ q}}\sum_{b\bmod\,q}f(b)\sum_{a\bmod\,q}\mathrm{e}\left(\dfrac{a(n\pm b)}{q}\right)\\
&=&f(\mp n),
\end{eqnarray*}
by the orthogonality relation for additive characters (or simply using geometric series). This proves 
\begin{equation}\label{ftf}
\left(f^{\pm}\right)^+=f_\mp=\delta f.
\end{equation}
Inserting this into \eqref{neu} leads to
\begin{equation}\label{ldel}
\Delta(1-s)\Delta(s)=\delta.
\end{equation}
Thus using \eqref{two}, we get
\begin{eqnarray}\label{i3}
{\mathcal I}_3 &=& -\dfrac{1}{ 2\pi i}\int_{-\epsilon+it_a}^{-\epsilon+iT} {\dfrac{\Delta'}{\Delta}}(s)\left(1+\dfrac{a}{\Delta(s)}+\sum_{j\geq 2}\left(\dfrac{a}{\Delta(s)}\right)^j\right)L(s;f)\d s \nonumber\\
&=&\dfrac{1}{ 2\pi i}\int_{1+\epsilon-it_a}^{1+\epsilon-iT} \dfrac{\Delta'}{\Delta}(1-s) \nonumber \\
&&\qquad \times \left(1+\delta a\Delta(s)+\sum_{j\geq 2}\left(\delta a\Delta(s)\right)^j\right)L(1-s;f)\d s\nonumber \\
&=:&\sum_{1\leq \ell\leq 3}\mathcal{J}_\ell.
\end{eqnarray}

We estimate at first $\mathcal{J}_1$ by computing its conjugate $\overline{\mathcal{J}_1}$. The functional equation \eqref{feq2}, relations \eqref{delta1} and \eqref{refle}, estimate \eqref{log_der}, and the Dirichlet series representation of $L\left(s;\overline{f^+}\right)$ in the half-plane $\sigma>1$ yield
\begin{eqnarray}\label{grobi}
\overline{\mathcal J}_1 &=& -\dfrac{1}{ 2\pi i}\overline{\int_{t_a}^{T}{\dfrac{\Delta'}{\Delta}(-\epsilon+it)}{\Delta(-\epsilon+it)}{L\left(1+\epsilon-it;{f}^+\right)}(-i)\d t}\nonumber\\
&=&\dfrac{1}{2\pi i}{\int_{t_a}^{T}-\overline{\dfrac{\Delta'}{\Delta}(-\epsilon+it)}\,\delta{\Delta(-\epsilon-it)}{L\left(1+\epsilon+it;\overline{{f}^+}\right)}i\d t}\nonumber\\
&=& \delta\int_{t_a}^{T}\left(\log\dfrac{q\tau}{2\pi}+O\left(\dfrac{1}{\tau}\right)\right)\d\left(\frac{1}{2\pi i}\int_{1+\epsilon+i}^{1+\epsilon+i\tau}\Delta(1-s)L\left(s;\overline{f^+}\right)\d s\right).
\end{eqnarray}
We evaluate the inner integral by applying Gonek's lemma: {\it Suppose that $\sum_{n=1}^{\infty}a(n)n^{-s}$ converges for $\sigma>1$ where $a(n)\ll n^\epsilon$ for any $\epsilon >0$. 
Let $\omega=\pm 1$ and $\epsilon>0$. 
Then we have}
\begin{eqnarray*}
&&\frac{1}{2\pi i}\int_{1+\epsilon+i}^{1+\epsilon+i\tau}\left(\frac{q}{2\pi}\right)^s\Gamma(s)\exp\left(\omega \frac{\pi i s}{2}\right)\sum_{n=1}^{\infty}\frac{a(n)}{n^s}\d s \\
&&\qquad= \begin{cases}
\displaystyle \sum_{n\leq \frac{\tau q}{2\pi}}a(n)\exp(-2\pi i\frac{n}{q})+O_\epsilon\left(\tau^{1/2+\epsilon}\right),&\mbox{ if } \quad \omega = -1, \\[5mm]
O_\epsilon(1), & \mbox{ if }\quad \omega = +1.
\end{cases}
\end{eqnarray*}
A proof of this lemma follows along the lines of \cite[Lemma 3]{steudi1}. A slightly weaker statement was proven originally by Gonek \cite[Lemma 5]{gonek} (and there exists a preversion of it in \cite[\S 7.4]{tit}).
Since
$$
\Delta(1-s)=\left(\dfrac{q}{2\pi}\right)^s\dfrac{\Gamma(s)}{\sqrt{q}}\left(\delta \mathrm{e}\left(\frac{s}{4}\right)+\mathrm{e}\left(-\frac{s}{4}\right)\right)
$$
by the definition (\ref{delta}) of $\Delta$, this gives here
$$
\frac{1}{2\pi i}\int_{1+\epsilon+i}^{1+\epsilon+i\tau}\Delta(1-s)L\left(s;\overline{f^+}\right)\d s
=\dfrac{1}{ \sqrt{q}}\sum_{n\leq \frac{\tau q}{ 2\pi}}\overline{f^+}(n)\,\mathrm{e}\left(-\dfrac{n}{q}\right)+O\left(\tau^{1/2+\epsilon}\right);
$$
we can simplify the right hand side: taking into account the $q$-periodicity, we find
\begin{eqnarray*}
\sum_{n\leq \frac{\tau q}{ 2\pi}}\overline{f^+}(n)\,\mathrm{e}\left(-\dfrac{n}{q}\right)&=&\sum_{a\bmod\,q}\overline{f^+}(a)\,\mathrm{e}\left(-\dfrac{a}{q}\right)\mathop{\sum_{n\leq \frac{\tau q}{ 2\pi}}}_{ n\equiv a\bmod\,q}1\\
&=&\overline{\sum_{a\bmod\,q}{f}^+(a)\,\mathrm{e}\left(\dfrac{a}{q}\right)}\left(\dfrac{\tau}{ 2\pi}+O(1)\right)\\
&=&\sqrt{q}\,\overline{(f^+)^+}(1)\dfrac{\tau}{ 2\pi}+O(1) \\
&=& \sqrt{q}\delta\,\overline{f}(1)\dfrac{\tau}{ 2\pi}+O(1).
\end{eqnarray*}
This leads via (\ref{grobi}) to
\begin{equation}\label{j1}
\mathcal{J}_1 = f(1)\frac{T}{2\pi}\log\frac{qT}{2\pi e}+O\left(T^{1/2+\epsilon}\right).
\end{equation}

Next we consider 
$$
{\mathcal J}_2 = \dfrac{a}{ 2\pi\delta i}\int_{1+\epsilon-it_a}^{1+\epsilon-iT}\frac{\Delta'}{ \Delta}(1-s)\Delta(s)\Delta(1-s)L\left(s;f^+\right)\d s.
$$
In view of (\ref{ldel}) this equals $a/(2\pi)$ times the conjugate of 
\begin{eqnarray*}
\lefteqn{-\int_{t_a}^{T}\overline{\dfrac{\Delta'}{\Delta}(-\epsilon+it)L(1+\epsilon-it;{f^+})}\d t}\\&=&\int_{t_a}^T\left(\log\dfrac{qt}{ 2\pi}+O\left(\dfrac{1}{t}\right)\right)\sum_{n\geq 1}\dfrac{\overline{f^+}(n)}{n^{1+\epsilon+it}}\d t\\
&=&\overline{f^+}(1)\int_{t_a}^{T}\log\frac{qt}{2\pi}\d t+\sum_{n\geq2}\dfrac{\overline{f^+}(n)}{n^{1+\epsilon}} J_n(T)+O(1),
\end{eqnarray*}
where we have used the absolute convergence in $\sigma>1$ and $J_n(T)$ is for $n\geq 2$ given by 
$$
J_n(T):=\int_{t_a}^T\log\dfrac{qt}{ 2\pi}\exp(-it\log n)\d t.
$$
To bound this integral we shall use the first derivative test: {\it Given real functions $F$ and $G$ on $[a,b]$ such that $G(t)/F'(t)$ is monotonic and $F'(t)/G(t)\geq M>0$ or $F'(t)/G(t)\leq -M<0$, then}
$$
\int_a^bG(t)\exp(iF(t))\d t\ll \dfrac{1}{ M}.
$$
This is essentially a classical lemma from \cite[\S 4.3]{tit}. 
This leads to the estimate $J_n(T)\ll \frac{\log T}{\log n}$ for $n\geq 2$. 
Hence, the contribution of the tail of the Dirichlet series is negligible and we get
\begin{equation}\label{j2}
{\mathcal J}_2=af^+(1)\dfrac{T}{ 2\pi}\log\dfrac{qT}{ 2\pi e}+O(\log T). 
\end{equation}

Finally, we have to consider the third integral over the tail of the geometric series. 
For this aim we apply \eqref{pow}, \eqref{log_der} and \eqref{oio} and get 
$$
\mathcal{J}_3 \ll \int_{t_a}^{T} \log t\sum_{j\geq 2}t^{-j(1/2+\epsilon)}t^{1/2+\epsilon}\d t \ll T^{1/2+\epsilon}.
$$
This together with \eqref{j1} and \eqref{j2} substituted into \eqref{i3}, in combination with \eqref{i1} shows that the vertical integrals contribute
\begin{eqnarray}\label{verti}
\mathcal{I}_1+\mathcal{I}_3&=&\left(f(1)+af^+(1)\right)\dfrac{T}{ 2\pi}\log\dfrac{qT}{ 2\pi e}+O(T^{1/2+\epsilon}),\end{eqnarray}
which is already the main term. 
It remains to consider the horizontal integrals in (\ref{i1234}). 
Although they can be treated the same way as in \cite{steusur} we sketch the details.

The integral $\mathcal{I}_2$ in \eqref{i1234} can be rewritten with the aid of the truncated partial fraction decomposition \eqref{aprox} as
$$
{\mathcal I}_2= \dfrac{1}{ 2\pi i}\int_{2+iT}^{-\epsilon+iT}\left(\sum_{\vert T-\gamma_a\vert\leq 1}\frac1{s-\delta_a}+O(\log T)\right)L(s;f)\d s. 
$$
Taking into account (\ref{condi}) we notice that $1/\vert s-\delta_a\vert\ll \log T$. 
Hence, utilizing \eqref{app1}, \eqref{ioi} and \eqref{oio}, we can show
\begin{eqnarray*}
{\mathcal I}_2&\ll &(\log T)^2 \left\{\int_{-\epsilon}^{0}+\int_{0}^{1}+\int_1^{2}\right\}\vert L\left(\sigma+iT;f\right)\vert\d \sigma\\
&\ll_\epsilon & (\log T)^2\left\{T^{1/2+\epsilon}+T^{1/2+\epsilon}+1\right\}\\
&\ll_\epsilon& T^{1/2+\epsilon}, 
\end{eqnarray*}
where $\epsilon$ at different places may take different values. This is the bound for the horizontal integrals. In combination with \eqref{verti} we arrive via \eqref{i1234} at the asymptotic formula of the theorem. Taking into account \eqref{ioi} and \eqref{oio} we may replace the chosen $T$ (with respect to \eqref{condi}) with a general $T\geq t_a$ at the expense of an error of order $T^{1/2+\epsilon}$ (as follows from \eqref{ioi}). This concludes the proof of Theorem \ref{mean}.

\section{Proof of Theorem \ref{uniform}}

We first note that \eqref{pow} implies
\begin{equation} \label{delta_asymp}
\Delta(\sigma+it) \asymp \left(\dfrac{qt}{2\pi}\right)^{1/2-\sigma},\quad t\geq t_0>0,
\end{equation}
while \eqref{amod} yields
\begin{equation}\label{rp}
\beta_a=\dfrac{1}{2}-\dfrac{\log|a|}{\log\frac{q\gamma_a}{2\pi}}+O\left(\dfrac{1}{\gamma_a\log\frac{q\gamma_a}{2\pi}}\right)
\end{equation}
for any $\gamma_a\geq T_0:=\max\left\{1,4\pi/q\right\}$.

Let $T\geq T_0$, $x>1$ and
$$
\alpha=\alpha(T) := \dfrac12+\dfrac{c}{\log{\frac{qT}{2\pi}}}
$$
where $c>0$ is a sufficiently large constant depending on $a$ satisfying
\begin{equation}\label{dist}
\left|\beta_a-\dfrac{1}{2}\right|\leq\dfrac{c}{2\log \frac{qT}{2\pi}}
\end{equation}
for any $a$-point with $\gamma_a\geq T\geq T_0$ as follows from \eqref{rp}, and
\begin{equation}\label{ineq}
\begin{array}{ccccccccc}
&|\Delta(\alpha+it)|&\leq &K\left(\dfrac{qt}{2\pi}\right)^{-c/\log \frac{qT}{2\pi}}&\leq &Ke^{-c}&<&\dfrac{|a|}{30}\\
\\
|\Delta(\alpha+it)|^{-1}= &|\Delta(1-\alpha+it)|&\geq& L\left(\dfrac{qt}{2\pi}\right)^{c/\log \frac{qT}{2\pi}}&\geq &Le^{c}&>&30|a|
\end{array}
\end{equation}
for any $t\geq T\geq T_0$, where $K,L$ are absolute constants coming from \eqref{delta_asymp} and we can assume without loss of generality that $K>1$ is sufficiently large and $L=1/K$.
After choosing a suitable $K$, we then take $c$ sufficiently large. These show that $\Delta(\alpha+it) \neq a$, $\Delta(1-\alpha+it) \neq a$ when $t\ge T$.

We know from \eqref{na} that for any $T_0\leq T<T+1\leq T'\leq2T$, we can find $T_1\in\left[T,T+1/2\right)$ and $T_2\in\left(T'-1/2,T'\right]$ such that
\begin{align}\label{T_i}
\begin{array}{ccc}\min\limits_{\ell=1,2}\min\limits_{\delta_a}\left|T_\ell-\gamma_a\right|\gg\dfrac{1}{\log T}\end{array}.
\end{align}

If $\mathcal{C}$ is the positively oriented rectangular contour with vertices $\alpha+iT_1$, $\alpha+iT_2$, $1-\alpha+iT_2$, $1-\alpha+iT_1$, then \eqref{na}, \eqref{dist} and the calculus of residues yield
\begin{eqnarray}\label{anf}
 \sum_{\substack{T<\gamma_a<T'}} x^{\delta_a}&=& \sum_{\substack{T_1<\gamma_a<T_2}} x^{\delta_a}+\sum_{\substack{T<\gamma_a<T_1}}x^{\delta_a}+\sum_{\substack{T_2<\gamma_a<T'}} x^{\delta_a}\nonumber\\
 &=&\frac{1}{ 2\pi i} \int_{\mathcal{C}} \frac{\Delta'(s)}{\Delta(s)-a} x^s\, \d s+O\left(x^\alpha\log T\right).
\end{eqnarray}

We break the integral
$$ \frac{1}{ 2\pi i} \int_{\mathcal{C}} \frac{\Delta'(s)}{\Delta(s)-a} x^s\, \d s $$
down into
$$
\frac{1}{2\pi i} \left(\int_{\alpha+iT_1}^{\alpha+iT_2} + \int_{\alpha+iT_2}^{1-\alpha+iT_2} + \int_{1-\alpha+iT_2}^{1-\alpha+iT_1} + \int_{1-\alpha+iT_1}^{\alpha+iT_1} \right) \frac{\Delta'(s)}{\Delta(s)-a} x^s \,\d{s} =: \sum_{1\leq j\leq 4}I_j.
$$

In view of \eqref{aprox} and \eqref{T_i}, for $\ell= 1, 2$ we have
\begin{eqnarray*}
\int_{1-\alpha+iT_\ell}^{\alpha+iT_\ell} \frac{\Delta'(s)}{\Delta(s)-a} x^s \,\d{s}
&=& \int_{1-\alpha+iT_\ell}^{\alpha+iT_\ell} \left(\sum_{\vert t-\gamma_a\vert\leq 1}\frac1{s-\delta_a} + O(\log{t})\right) x^s \,\d{s} \\
&=& \int_{1-\alpha+iT_\ell}^{\alpha+iT_\ell} \sum_{\vert t-\gamma_a\vert\leq 1} \frac1{s-\delta_a} x^s \,\d{s}
+ O\left(x^\alpha\right),
\end{eqnarray*}
since
$$
\int_{1-\alpha+iT_\ell}^{\alpha+iT_\ell} O(\log{t}) x^s \,\d{s}
\ll \log T_\ell \int_{1-\alpha}^\alpha x^\sigma\, \d\sigma
\ll x^\alpha (2\alpha-1) \log T \ll x^\alpha.
$$
Meanwhile,
\begin{eqnarray*}
\int_{1-\alpha+iT_\ell}^{\alpha+iT_\ell} \sum_{\vert t-\gamma_a\vert\leq 1} \frac1{s-\delta_a} x^s \,\d{s}
&\ll& \sum_{\vert T_\ell-\gamma_a\vert\leq 1} \int_{1-\alpha}^{\alpha} \frac{x^\sigma}{\left|\alpha-\beta_a+i(T_\ell-\gamma_a)\right|} \,\d{\sigma} \\
&\ll& x^\alpha (2\alpha-1) \log T \sum_{\vert T_1-\gamma_a\vert\leq 1} 1 \\
&\ll& x^\alpha\log T.
\end{eqnarray*}
Therefore,
\begin{equation} \label{hor}
I_2,I_4 \ll x^\alpha\log{T}.
\end{equation}

We now estimate the vertical integrals $I_1$ and $I_3$. For $\sigma=\alpha$, we use \eqref{bibo2} and thus
\begin{equation} \label{delta_i1_rewrite}
\frac{\Delta'(s)}{\Delta(s)-a} = \frac{\Delta'(s)}{ -a} \cdot \frac{1}{1 - \Delta(s)/a}
= \frac{\Delta'(s)}{ -a} \left( 1 + \sum_{j\ge1} \left(\frac{\Delta(s)}{a}\right)^j \right);
\end{equation}
here the first inequality in (\ref{ineq}) allows us to expand the second factor into a geometric series. Applying this, we find
\begin{equation}\label{I_1ref}
I_1 = -\frac{x^{\alpha}}{2\pi a} \int_{T_1}^{T_2} \Delta'(\alpha+it) \left( \sum_{0\leq j< m} \left(\frac{\Delta(\alpha+it)}{a}\right)^j + \sum_{j\ge m} \left(\frac{\Delta(\alpha+it)}{a}\right)^j \right) x^{it} \,\d{t}.
\end{equation}

Here we pick $m$ depending on $T$, large enough to bound the last term in the integral trivially. By \eqref{delta_asymp} and again the first inequality in \eqref{ineq},
\begin{equation} \label{delta_i1_asymp}
\sum_{j\ge m} \left(\frac{\Delta(\alpha+it)}{a}\right)^j
= \frac{ \left(\frac{\Delta(\alpha+it)}{a} \right)^m}{1-\frac{\Delta(\alpha+it)}{a}}
\ll\dfrac{30}{29}\left(\dfrac{1}{30}\right)^m.
\end{equation}
Take
\begin{equation} \label{m_of_T}
m := \left\lfloor\dfrac{2\log{\frac{qT}{4\pi}}}{\log 30}\right\rfloor.
\end{equation}
It follows then by \eqref{delta_i1_asymp} that
\begin{eqnarray*}
I_1 &=& -\frac{x^{\alpha}}{2\pi a} \int_{T_1}^{T_2} \Delta'(\alpha+it) \sum_{0\leq j< m} \left(\frac{\Delta(\alpha+it)}{a}\right)^jx^{it} \,\d{t} \\
&&+~ O\left(x^\alpha\int_{T_1}^{T_2}\left|\Delta'(\alpha+it)\right|\left(\dfrac{4\pi}{qT}\right)^{2} \d{t}\right).
\end{eqnarray*}
The first term of the integrand can be estimated using \eqref{pow} for which we obtain
\begin{eqnarray*}
\Delta'(\alpha+it) &=& \delta\left(\frac{qt}{ 2\pi}\right)^{1/2-(\alpha+it)}\exp\left({i\left(t+\frac{\pi}{ 4}\right)}\right) \left( -\log\frac{qt}{ 2\pi} + O\left(\frac{\log{qt}}{t}\right) \right)\\
&\ll& \log\frac{qt}{2\pi}.
\end{eqnarray*}
Here we substituted the value of $\alpha$ and used the first inequality in \eqref{ineq}. 
It then follows that we can discard the last term in $I_1$ as
$$
I_1 = -\frac{x^{\alpha}}{2\pi a} \int_{T_1}^{T_2} \Delta'(\alpha+it) \sum_{0\le j< m} \left(\frac{\Delta(\alpha+it)}{a}\right)^j x^{it} \,\d{t} + O\left(x^{\alpha}\right).
$$
Since
$$
\frac\d{\d{t}} \Delta(\alpha+it)^{j+1} = i(j+1)\Delta'(\alpha+it)\Delta(\alpha+it)^j,
$$
we can rewrite $I_1$ as
\begin{eqnarray}\label{I_1}
I_1 &=& -\frac{x^{\alpha}}{2\pi i}\sum_{1\le j\le m} \dfrac{1}{ja^j} \int_{T_1}^{T_2} \left(\Delta(\alpha+it)^{j}\right)' x^{it} \,\d{t} + O\left(x^{\alpha}\right)\nonumber\\
&=:& -\frac{x^{\alpha}}{2\pi i}\sum_{1\le j\le m} \dfrac{1}{ja^j} I_{1j}+ O\left(x^{\alpha}\right).
\end{eqnarray}

We estimate $I_{1j}$ for $1\le j\le m$. Integrating by parts, we obtain with the aid of \eqref{pow} and the first inequality in \eqref{ineq}
\begin{eqnarray}
I_{1j}&=&\Delta(\alpha+it)^{j} x^{it}\Big|_{T_1}^{T_2}-i\log x\int_{T_1}^{T_2}\Delta(\alpha+it)^{j}x^{it}\d{t}\nonumber\\
&\ll&\log x\left|\int_{T_1}^{T_2}\delta^j\left(\frac{qt}{ 2\pi}\right)^{(1/2-\alpha-it)j}\exp\left(ij\left(t+\dfrac{\pi}{ 4}\right)\right)\left(1+O\left(t^{-1}\right)\right)^jx^{it}\d{t}\right| \nonumber\\
&&+ \left(\dfrac{|a|}{30}\right)^j \nonumber.
\end{eqnarray}
We see that if $|O\left(t^{-1}\right)|<D/t$ for some $D>1$, then
$$ 
\left| \left(1+O\left(t^{-1}\right)\right)^j-1 \right|
= \left| \sum\limits_{k=1}^j\binom{j}{k}\left(O\left(t^{-1}\right)\right)^k \right|
\leq \frac{D^j}t \sum\limits_{k=1}^j\binom{j}{k} \leq \frac{(2D)^j}t 
$$
or
$$ 
\left(1+O\left(t^{-1}\right)\right)^j=1+O\left((2D)^jt^{-1}\right). 
$$
In view of \eqref{ineq}, we then have for a sufficiently large constant $K$ that
\begin{eqnarray}\label{I_1j}
I_{1j} &\ll& \log x\left|\int_{T_1}^{T_2}\left(\frac{qt}{ 2\pi}\right)^{\left(-c/\log \frac{qT}{2\pi}-it\right)j} \exp\left(ijt\right) \left(x^{1/j}\right)^{ijt} \d{t} \right| \nonumber\\
&&+\, \log x\int_{T_1}^{T_2}\left(2D\left(\frac{qt}{ 2\pi}\right)^{-c/\log \frac{qT}{2\pi}}\right)^jt^{-1}\d{t}
+ \left(\dfrac{|a|}{30}\right)^j \nonumber\\
&\ll& q^{-jc/\log \frac{qT}{2\pi}} \log{x} \left| \int_{T_1}^{T_2} \left(\frac{t}{ 2\pi}\right)^{-jc/\log \frac{qT}{2\pi}}
\exp\left[-ijt\log\left(\dfrac{qt}{2\pi x^{1/j}e}\right) \right]\d{t}\right| \nonumber\\
&&+ \left(2D\dfrac{|a|}{30K}\right)^j \log x \int_{T_1}^{T_2}t^{-1}\d{t}+\left(\dfrac{|a|}{30}\right)^j\nonumber\\
&\ll&\dfrac{\log x}{j}\left(\dfrac{q}{j}\right)^{-jc/\log \frac{qT}{2\pi}}\left|\mathcal{J}\right| + (1+\log x)\left(\dfrac{|a|}{30}\right)^j,
\end{eqnarray}
where
$$
\mathcal{J}:=\int_{jT_1}^{jT_2}\exp\left[-it\log\left(\dfrac{qt}{2\pi jx^{1/j}e} \right)\right]\left(\frac{t}{ 2\pi}\right)^{1/2-jc/\log \frac{qT}{2\pi}-1/2}\d{t}.
$$
To estimate $\mathcal{J}$ we employ another form of Gonek's lemma (see \cite[Lemma 2]{gonek}):
{\it For large $A$ and $A<r\leq B\leq 2A$
\begin{eqnarray*}
\lefteqn{\int_A^B\exp\left[it\log\left(\dfrac{t}{re}\right)\right]\left(\dfrac{t}{2\pi}\right)^{\mathfrak{a}-1/2}\d{t}}\\
&=& (2\pi)^{1-\mathfrak{a}}r^\mathfrak{a}\exp\left(-i\Big(r-{\pi\over 4}\Big)\right)\mathbbm{1}_{[A,B)}(r)+E(r,A,B),
\end{eqnarray*}
where $\mathfrak{a}$ is a fixed real number and
$$
E(r,A,B) \ll A^{\mathfrak{a}-1/2}+\dfrac{A^{\mathfrak{a}+1/2}}{|A-r|+A^{1/2}}+\dfrac{B^{\mathfrak{a}+1/2}}{|B-r|+B^{1/2}}.
$$}
It is easily seen that that this holds uniformly in $\mathfrak{a}$ from a bounded interval (the implicit constants depend of course on the limits of this interval). We then apply this, by taking complex conjugate, to $\mathcal{J}$ with $\mathfrak{a}=1/2-jc/\log \frac{qT}{2\pi}$, $r=2\pi jx^{1/j}/q$, $A=jT_1$ and $B=jT_2$:
\begin{eqnarray*}
\mathcal{J}&\ll&(2\pi)^{1/2+jc/\log \frac{qT}{2\pi}}\left(\dfrac{2\pi jx^{1/j}}{q}\right)^{1/2-jc/\log \frac{qT}{2\pi}}\mathbbm{1}_{\left[jT_1,jT_2\right)}\left(\dfrac{2\pi jx^{1/j}}{q}\right) \nonumber\\
&&+\, \left(jT_1\right)^{-jc/\log \frac{qT}{2\pi}}
+ \dfrac{\left(jT_1\right)^{1-jc/\log \frac{qT}{2\pi}}}{\left|jT_1-\frac{2\pi jx^{1/j}}{q}\right|+\sqrt{jT_1}}
+ \dfrac{\left(jT_2\right)^{1-jc/\log \frac{qT}{2\pi}}}{\left|jT_2-\frac{2\pi jx^{1/j}}{q}\right|+\sqrt{jT_2}} \nonumber\\
&\ll& \left(\dfrac{j}{q}\right)^{1/2-jc/\log \frac{qT}{2\pi}}x^{1/(2j)-c/\log \frac{qT}{2\pi}}
\mathbbm{1}_{\left(T/2,{5}T/2\right)}\left(\dfrac{2\pi x^{1/j}}{q}\right) \\
&&+~ (jT)^{-jc/\log \frac{qT}{2\pi}}E_0(x,j,T), \nonumber
\end{eqnarray*}
where
\begin{eqnarray}\label{error}
E_0(x,j,T)
&=&\left\{\begin{array}{ll}
O\left(\sqrt{jT}\right),
&\text{if } \dfrac{2\pi x^{1/j}}{q}\in\left(\dfrac{T}{2},\dfrac{{5}T}{2}\right), \\[3mm]
O(1), &\text{otherwise}.
\end{array}\right.
\end{eqnarray}
Recall that $T\asymp T_1\asymp T_2$ and $j\ll\log T$. Therefore, $E_0$ does not restrictively depend on $T_1$ and $T_2$.

Applying this to \eqref{I_1j} we have
\begin{eqnarray*}
I_{1j} &\ll& \dfrac{\log x}{\sqrt{jq}}x^{1/(2j)-c/\log \frac{qT}{2\pi}}\mathbbm{1}_{\left(T/2,{5}T/2\right)}\left(\dfrac{2\pi x^{1/j}}{q}\right) \\
&&+~ \dfrac{\log x}{j}\left(\left(qT\right)^{-c/\log \frac{qT}{2\pi}}\right)^jE_0(x,j,T)
+ (1+\log x) \left(\dfrac{|a|}{30}\right)^j.
\end{eqnarray*}
Hence, in view of \eqref{I_1} and \eqref{ineq}, we obtain
\begin{eqnarray}\label{sum1}
I_1&\ll& x^{1/2}\log x\sum\limits_{1\leq j\le m} \dfrac{x^{1/(2j)}}{j^{3/2}|a|^j}\mathbbm{1}_{\left(T/2,{5}T/2\right)}\left(\dfrac{2\pi x^{1/j}}{q}\right) \nonumber\\
&&+~ x^\alpha\log x\sum\limits_{1\leq j\le m}\dfrac{1}{j^2|a|^j}\left(\left(2\pi\right)^{-c/\log \frac{qT}{2\pi}}\dfrac{1}{K}\dfrac{|a|}{30}\right)^jE_0(x,j,T) \nonumber\\
&&+~ x^\alpha(1+\log x)\sum\limits_{1\leq j\le m}\dfrac{1}{j30^j} + x^\alpha.
\end{eqnarray}
Observe that the intervals
\begin{eqnarray*}
\left(\left(\dfrac{qT}{4\pi}\right)^j,\left(\dfrac{{5}qT}{4\pi}\right)^j\right), &1\le j\le m,
\end{eqnarray*} are pairwise disjoint whenever
$$ j<\dfrac{1}{\log{5}}\log\left(\dfrac{qT}{4\pi}\right). $$
Comparing this with \eqref{m_of_T}, we see that
$$ m \leq \dfrac{2}{\log{30}}\log\left(\dfrac{qT}{4\pi}\right)< \dfrac{1}{\log{5}}\log\left(\dfrac{qT}{4\pi}\right) $$
for all $T\geq T_0$.
By this construction, the first sum in \eqref{sum1} can have at most only one term, namely,
$$\dfrac{x^{1/(2j_x)}}{j_x^{3/2}|a|^{j_x}}$$
for the unique, if any, $j_x\le m$ such that $2\pi x^{1/j_x}/q\in(T/2,{5}T/2)$.
Similar reasoning also applies to $E_0$ of the second sum in \eqref{sum1} and it follows in view of \eqref{error} that
$$ \sum\limits_{1\leq j\le m} \dfrac{1}{j^2 30^j} E_0(x,j,T)
\ll \mathop{\sum\limits_{1\leq j\le m}}_{j\neq j_x} \dfrac{1}{j^230^j}
+ \dfrac{T^{1/2}}{j_x^{3/2}30^{j_x}}\ll1+\dfrac{T^{1/2}}{j_x^{3/2}30^{j_x}}. $$
The last sum in \eqref{sum1} is trivially $O(1)$.
Collecting these estimates, we conclude that
\begin{equation}\label{I_1f}
I_1\ll x^{\alpha}(1+\log x)+E_1(x,a,T),
\end{equation}
where $E_1(x,a,T)$ is defined to be equal to
\begin{equation}\label{E''}
\dfrac{\log x}{j_x^{3/2}}\left(\dfrac{x^{1/2+1/(2j_x)}}{|a|^{j_x}}+\dfrac{x^{\alpha}T^{1/2}}{30^{j_x}}\right)
\end{equation}
if such $j_x$ exists, and zero otherwise.
Observe that if $T\geq T_0\geq4\pi/q$, then for any $x>1$ it is $E_1(x,a,T)\geq0$, while for $0<x<1$ we always have $E_1(x,a,T)=0$.

For $\sigma=1-\alpha$, we use \eqref{bibo} and the corresponding \eqref{two} to write
\begin{equation} \label{delta_i3_rewrite}
\frac{\Delta'(s)}{\Delta(s)-a} = \frac{\Delta'}{\Delta}(s) \cdot \frac{1}{1 - a/\Delta(s)}
=\frac {\Delta'}{\Delta}(s) \left( 1 + \sum_{j\ge1} \left(\frac{a}{\Delta(s)}\right)^j \right);
\end{equation}
here the second inequality in (\ref{ineq}) allows us to expand the second factor into a geometric series.
Thus the left vertical integral $I_3$ can be decomposed as follows
\begin{eqnarray*}
I_3 &= &\frac{1 }{ 2\pi i} \int_{1-\alpha+iT_2}^{1-\alpha+iT_1} \frac{\Delta'(s)}{\Delta(s)-a} x^s \,\d{s} \\
&= &\dfrac{1}{2\pi i}\int_{1-\alpha+iT_2}^{1-\alpha+iT_1}\dfrac{\Delta'}{\Delta}(s)x^s \d{s}
+ \dfrac{1}{2\pi i}\int_{1-\alpha+iT_2}^{1-\alpha+iT_1}\dfrac{\Delta'}{\Delta}(s)x^s\sum_{j\ge1} \left(\frac{a}{\Delta(s)}\right)^j \d{s} \\
&=:& I_{31}+I_{32}.
\end{eqnarray*}
Integrating by parts, we obtain in view of \eqref{log_der}
\begin{eqnarray}\label{I_31f}
I_{31} &=& \dfrac{x^{1-\alpha}}{2\pi} \int_{T_2}^{T_1} \left(-\log\dfrac{qt}{2\pi} + O\left(\dfrac{1}{t}\right)\right)x^{it}\d{t} \nonumber\\
&=& \left.\dfrac{x^{1-\alpha+it}}{2\pi i\log x} \left(\log\dfrac{qt}{2\pi} + O\left(\dfrac{1}{t}\right)\right) \right|_{T_1}^{T_2}
- \dfrac{x^{1-\alpha}}{2\pi i\log x}\int_{T_1}^{T_2}O\left(\dfrac{1}{t}\right)\d{t} \nonumber\\
&\ll&\dfrac{x^{{1-}\alpha}\log T}{\log x}.
\end{eqnarray}
In the case of $I_{32}$, it follows from \eqref{ldel} that
\begin{eqnarray*}
I_{32}&=&-\dfrac{1}{2\pi i}\int_{\alpha-iT_2}^{\alpha-iT_1}\dfrac{\Delta'}{\Delta}(1-s)x^{1-s}\sum_{j\ge1} \left(\frac{a}{\Delta(1-s)}\right)^j \d{s}\nonumber\\
&=&\dfrac{1}{2\pi i}\int_{\alpha-iT_1}^{\alpha-iT_2}\dfrac{\Delta'}{\Delta}(s)x^{1-s}\sum_{j\ge1} \left(\frac{a\Delta(s)}{\delta}\right)^j \d{s}\nonumber\\
&=&-\dfrac{x^{1-\alpha}}{2\pi}\int_{T_1}^{T_2}\dfrac{\Delta'}{\Delta}(\alpha-it)\sum_{j\ge1} \left(\frac{a\Delta(\alpha-it)}{\delta}\right)^j x^{it}\d{t}\nonumber\\
\end{eqnarray*}
or from \eqref{refle}
\begin{eqnarray}
\overline{I_{32}}&=&-\dfrac{x^{1-\alpha}}{2\pi}\int_{T_1}^{T_2}\dfrac{\overline{\Delta'(\alpha-it)}}{\overline{\Delta(\alpha-it)}}\sum_{j\ge1} \left(\frac{\overline{a}\overline{\Delta(\alpha-it)}}{\delta}\right)^jx^{-it} \d{t}\nonumber\\
&=&-\dfrac{\overline{a}x^{1-\alpha}}{2\pi}\int_{T_1}^{T_2}\Delta'(\alpha+it)\sum\limits_{j\geq0}\left(\overline{a}\Delta(\alpha+it)\right)^j\left(\dfrac{1}{x}\right)^{it}\d{t}.
\end{eqnarray}

To estimate now $\overline{I_{32}}$ we proceed exactly as the estimation of $I_1$ in \eqref{I_1ref}, where we have $1/\overline{a}$ instead of $a$, $x^{1-\alpha}$ instead of $x^{\alpha}$ and $y:=1/x$ instead of $x$.
We can then derive as in \eqref{I_1ref}--\eqref{I_1f} that
\begin{equation}\label{I_3f}
I_{32}\ll x^{1-\alpha}(1+\log x) + E_2(x,a,T),
\end{equation}
where $E_2(x,a,T)$ is defined to be equal to
\begin{equation}\label{E'''}
\dfrac{\log x}{j_{y}^{3/2}} \left( x^{1/2-1/(2j_{y})}|a|^{j_{y}} + \dfrac{x^{1-\alpha}T^{1/2}}{30^{j_{y}}} \right)
\end{equation}
if such $j_y$ exists, and zero otherwise.
Observe that if $T\geq T_0\geq4\pi/q$, then for any $x>1$ it is always $E_2(x,a,T)=0$, while for $0<x<1$ it is $E_2(x,a,T)\leq0$.

Collecting the estimates in \eqref{anf}, \eqref{hor}, \eqref{I_1f}, \eqref{I_31f} and \eqref{I_3f} we obtain that
$$
\sum_{\substack{T<\gamma_a< T'}} x^{\delta_a}\ll x^\alpha\left(\log x+\log T+\dfrac{\log T}{\log x}\right) + E_1(x,a,T)
$$
for any $x>1$ and any $T_0\leq T<T+1\leq T'\leq2T$.

The case of $0<x<1$ follows from \eqref{refle}, \eqref{ldel} and the above estimate. Indeed, relations \eqref{refle} and \eqref{ldel} imply that a complex number $z$ is an $a$-point of $\Delta(s)$ (where $a\neq0$) if and only if the complex number $1-\overline{z}$ is a $b:=1/\overline{a}$-point of $\Delta$.
Thus if $0<x<1$, then
$$\dfrac{1}{x}\sum_{\substack{T<\gamma_a<T'}} x^{\delta_a}
= \sum_{\substack{T<\gamma_a<T'}}\left(\dfrac{1}{ x}\right)^{1-\delta_a}
= \overline{\sum_{\substack{T<\gamma_a<T'}}\left(\dfrac{1}{ x}\right)^{1-\overline{\delta_a}}}
= \overline{\sum_{\substack{T<\gamma_b<T'}}\left(\dfrac{1}{ x}\right)^{\delta_b}},$$
or
$$
\sum_{\substack{T<\gamma_a<T'}} x^{\delta_a}
\ll x^{1-\alpha}\left(-\log x+\log T-\dfrac{\log T}{\log x}\right)
+ x \left(E_1\left(\frac1x,b,T\right)\right).
$$
By symmetry, it follows easily from \eqref{E''} and \eqref{E'''} that
\begin{eqnarray*}
xE_1(1/x,b,T) = -E_2(x,a,T)\geq0.
\end{eqnarray*}
Hence,
\begin{eqnarray*}
\sum_{\substack{T<\gamma_a<T'}} x^{\delta_a} &\ll& \left(x^\alpha+x^{1-\alpha}\right)\left(|\log x|+\log T+\dfrac{\log T}{|\log x|}\right) \\
&&+\, \mathbbm{1}_{(1,+\infty)}(x) E_1(x,a,T)-\mathbbm{1}_{(0,1)}(x)E_2(x,a,T)
\end{eqnarray*}
for any $0<x\neq1$ and any $T_0\leq T<T+1\leq T'\leq2T$.

The last statement of the theorem follows by the above construction. If $x\neq1$ is such that $4\pi/(qT)\leq x\leq qT/4\pi$, then $E_1(x,a,T)=E_2(x,a,T)=0$ and 
\begin{eqnarray*}
\sum_{\substack{T<\gamma_a<T'}} x^{\delta_a} &\ll& x^{1/2}\left(x^{c/\log\frac{qT}{2\pi}}+x^{-c/\log\frac{qT}{2\pi}}\right)\left(|\log x|+\log T+\dfrac{\log T}{|\log x|}\right)\\
&\ll& x^{1/2}\left(1+\dfrac{1}{|\log x|}\right)\log T.
\end{eqnarray*}

\section{Proof of Theorem \ref{u.d.}}

The Riemann-von Mangoldt type formula \eqref{na} implies that
\begin{equation}\label{RvM}
\begin{array}{ccc}N_a(T;f)\sim\dfrac{T\log T}{2\pi}&\text{ and }&\gamma_{a}^{(n)}\sim\dfrac{2\pi n}{\log n}.\end{array}
\end{equation}
Therefore, there is a $K\in\mathbb{N}$ such that $2^{K-1}\leq\gamma_a^{{(2n)}}/\gamma_a^{{(n)}}\leq 2^K$ for all $n\in\mathbb{N}$.
If we now set $\delta_{a}^{{(n)}}:=\beta_{a}^{{(n)}}+i\gamma_{a}^{{(n)}}$ to be the $a$-point of $\Delta$ with ordinate $\gamma_{a}^{{(n)}}$, then for any integers $1\leq N< M\leq 2N$ and any positive number $x\neq1$, we have
\begin{equation}\label{fs}
\sum\limits_{N< n\leq M}x^{\delta_{a}^{{(n)}}}
=\sum\limits_{k\leq K}\sum\limits_{2^{k-1}\gamma_{a}^{{(N)}}<\gamma_a\leq\min\left\{ \gamma_{a}^{{(M)}},2^{k}\gamma_{a}^{{(N)}}\right\}}x^{\delta_a},
\end{equation}
where the inner sum on the right-hand side is zero when the set of indices of its summands is empty.
Since $4\pi/(qN)\leq x\leq qN/(4\pi)$, \eqref{RvM}, \eqref{fs} and Theorem 3 yield
\begin{align*}\sum\limits_{N< n\leq M}x^{\delta_{a}^{{(n)}}}&\ll x^{1/2}\left(1+\dfrac{1}{|\log x|}\right)\log N.
\end{align*}
We also have
$$
x^{-\beta_a^{{(M)}}}\ll x^{-1/2\pm c/\log\frac{qN}{2\pi}}\ll x^{-1/2}
$$
for any $M>N$. Using this and Abel's summation formula (summation by parts), we obtain
\begin{eqnarray}\label{rim}
\sum\limits_{N< n\leq 2N}x^{i\gamma_{a}^{{(n)}}}
&=&\sum\limits_{N< n\leq 2N}x^{-\beta_{a}^{{(n)}}}x^{\delta_{a}^{{(n)}}}\nonumber\\
&\ll &x^{-\beta_{a}^{{(2N)}}}\left|\sum\limits_{N< n\leq 2N}x^{\delta_{a}^{{(n)}}}\right|+\nonumber\\
&&+\max\limits_{N< M<2N}\left|\sum\limits_{N< n\leq M}x^{\delta_{a}^{{(n)}}}\right|\sum\limits_{N< M<2N}\left|x^{-\beta_{a}^{{(M+1)}}}-x^{-\beta_{a}^{{(M)}}}\right|\nonumber\\
&\ll & \left(1+\dfrac{1}{|\log x|}\right)\log N\left(1+\sum\limits_{N<M< 2N}\left|x^{\beta_{a}^{{(M+1)}}-\beta_{a}^{{(M)}}}-1\right|\right).
\end{eqnarray}
By (\ref{amod}) resp. (\ref{rp}),
$$\dfrac{\beta_{a}^{{(M+1)}}-\beta_{a}^{{(M)}}}{\log|a|}=\left[\left({\log\frac{q\gamma_a^{{(M)}}}{2\pi}}\right)^{-1}-\left({\log\frac{q\gamma_a^{{(M+1)}}}{2\pi}}\right)^{-1}+O\left(\dfrac{1}{\gamma_a^{{(M)}}\log{\gamma_a^{{(M)}}}}\right)\right]$$
and $4\pi/(qN)\leq x\leq qN/(4\pi)$. Hence, with the aid of \eqref{RvM} we can show
\begin{eqnarray*}
\lefteqn{\sum\limits_{N<M<2N}\left|x^{\beta_{a}^{{(M+1)}}-\beta_{a}^{{(M)}}}-1\right|}\nonumber\\
&\ll &|\log x|\sum\limits_{N<M<2N}\left[\left({\log\frac{q\gamma_a^{{(M)}}}{2\pi}}\right)^{-1}-\left({\log\frac{q\gamma_a^{{(M+1)}}}{2\pi}}\right)^{-1}+\dfrac{1}{\gamma_a^{{(M)}}\log\frac{q\gamma_a^{{(M)}}}{2\pi}}\right]\\
&\ll &|\log x|\left[\left({\log\frac{\gamma_a^{{(N+1)}}}{2\pi}}\right)^{-1}-\left({\log\frac{\gamma_a^{{(2N-1)}}}{2\pi}}\right)^{-1}+\dfrac{N}{\gamma_a^{{(N)}}\log\gamma_a^{{(N)}}}\right]\\
&\ll &|\log x|\left[\left(\log \gamma_a^{{(N)}}\right)^{-2}+1\right]\\
&\ll & |\log x|.
\end{eqnarray*}
This and inequality \eqref{rim} lead for fixed $x$ to
\begin{equation}\label{la1}
\sum\limits_{N< n\leq 2N}x^{i\gamma_{a}^{{(n)}}}\ll\left(\dfrac{1}{|\log x|}+|\log x|\right)\log N=o_x\left(\left(\log N\right)^2\right),
\end{equation}
which is the first statement of the theorem.

Let now $\alpha\neq0$ be a real number and $k\neq0$ be an integer. If $x=\exp(2\pi\alpha k)$, then for every $N\geq N_x:=4\pi x/q$ and $K_x:=\left\lfloor\frac{\log (N/N_x)}{\log 2}\right\rfloor$ we have that
\begin{eqnarray*}
\sum\limits_{n\leq N}\mathrm{e}\left(k\alpha\gamma_a^{{{(n)}}}\right)&=&\sum\limits_{n\leq\frac{N}{2^{K_x}}}\mathrm{e}\left(k\alpha\gamma_a^{{{(n)}}}\right)+\sum\limits_{k\leq K_x}\sum\limits_{\frac{N}{2^{k}}< n\leq\frac{N}{2^{k-1}}}x^{i\gamma_{a}^{{(n)}}}\\
&=&O(N_x)+\sum\limits_{k\leq K_x}o_x\left(\left(\log N\right)^2\right)\\
&=&o_x\left(\left(\log N\right)^3\right),
\end{eqnarray*}
as follows from \eqref{la1}.
Since $k\in\mathbb{Z}\setminus\lbrace0\rbrace$ can be chosen arbitrarily, the sequence $\alpha\gamma^{{(n)}}_a$, $n\in\mathbb{N}$, satisfies Weyl's Criterion (see \cite{weyl}):
{\it a sequence $(x_n)_{n\in\mathbb{N}}$ is uniformly distributed modulo 1 if
\begin{equation*}\sum\limits_{n\leq N}\mathrm{e}(kx_n)=o(N)
\end{equation*}
for any integer $k\neq0$}. This concludes the proof of the theorem.

\section{Proof of Theorem \ref{universal}}

If $\chi$ is a Dirichlet character mod $r$ and $Q>0$, we define the truncated and twisted Euler product
\begin{equation}\label{trE}
s\longmapsto L_Q(s;\chi):=\prod\limits_{p\leq Q}\left(1-\dfrac{\chi(p)}{p^{s}}\right)^{-1}
\end{equation}
for every $s\in\mathbb{C}$ with $\sigma>0$, where $p$ will denote from here on a prime number.

The proof of \eqref{joint} is similar to Reich's proof \cite{reich} for the discrete universality of $\zeta(s)$ and we will not repeat it here. Instead we employ the following theorem which highlights the necessary conditions a sequence $(x_n)_{n\in\mathbb{N}}$ has to meet in order to derive universality: {\it Let $\chi_1,\dots , \chi_J$ be pairwise non-equivalent Dirichlet characters. Let also $(x_n)_{n\in\mathbb{N}}$ be a sequence of real numbers such that the sequence of vectors
\begin{equation}\label{con1}
\begin{array}{cc}\left(x_n\dfrac{\log p}{2\pi}\right)_{p\in\mathcal{M}},&n\in\mathbb{N},
\end{array}
\end{equation}
is uniformly distributed modulo 1 for any finite set of primes $\mathcal{M}$, and
\begin{equation}\label{con2}\lim\limits_{Q\to\infty}\limsup\limits_{N\to\infty}\dfrac{1}{N+1}\sum\limits_{N\leq n\leq2N}\left|L\left(s+ix_n;\chi_j\right)-L_Q\left(s+ix_n;\chi_j\right)\right|^2=0,
\end{equation}
$j=1,\dots,J$, uniformly in compact subsets of the strip $1/2<\sigma<1$. Then for any compact subset with connected complement $K$ of this strip, any $g_1,\dots,g_J$ continuous non-vanishing functions on $K$ and analytic in its interior, any $z>0$ and $\xi_p$, $p\leq z$, real numbers, and any $\varepsilon>0$,}
$$
\liminf\limits_{N\to\infty}\dfrac{1}{N}\sharp\left\{1\leq n\leq N:\begin{array}{c}\max\limits_{1\leq j\leq J}\max\limits_{s\in K}\left|L\left(s+ix_n;\chi_j\right)-g_j(s)\right|<\varepsilon\\ \& \
\max\limits_{p\leq z}\left\|\gamma_a^{{(n)}}\dfrac{\log p}{2\pi}-\xi_p\right\|<\varepsilon
\end{array}\right\}>0.
$$
For the proof see \cite[Theorem 3.1]{sourm}.

The definition of uniform distribution of a multidimensional sequence (sequence of vectors) is analogous to the one-dimensional case (sequence of numbers) and so we omit the details here. 
However, we will use an equivalent statement of it (see \cite[Chapter I, Theorem 6.3]{kuinie}): {\it A sequence $\left(\underline{x}_n\right)_{n\in\mathbb{R}}$, of vectors from $\mathbb{R}^{\ell}$ (for some $\ell\in\mathbb{N}$) is uniformly distributed modulo 1 if, and only if, for every $\underline{h}\in\mathbb{Z}^\ell\setminus\lbrace\underline{0}\rbrace$, the sequence $\left\langle\underline{h},\underline{x}_n\right\rangle$, $n\in\mathbb{N}$, is uniformly distributed modulo $1$, where $\langle\cdot,\cdot\rangle$ is the standard inner product.}

\noindent Therefore, if $\mathcal{M}$ is a finite set of primes, then the sequence
 \begin{eqnarray*}
 \left(\gamma_a^{{(n)}}\dfrac{\log p}{2\pi}\right)_{p\in \mathcal{M}},&n\in\mathbb{N},
 \end{eqnarray*}
 is uniformly distributed modulo $1$ if and only if the sequence
 \begin{eqnarray*}
 \gamma_a^{{(n)}}\sum\limits_{p\in \mathcal{M}}h_p\dfrac{\log p}{2\pi},&n\in\mathbb{N},
 \end{eqnarray*}
 is uniformly distributed modulo $1$ for any
 $(h_p)_{p\in\mathcal{M}}\in\mathbb{Z}^{\sharp\mathcal{M}}\setminus\lbrace\underline{0}\rbrace$. 
But this follows immediately from Theorem \ref{u.d.} and the unique factorization of integers into primes.
Thus, $\gamma_{a}^{{(n)}}$, $n\in\mathbb{N}$, satisfies condition \eqref{con1}.
 
To prove that $\gamma_{a}^{{(n)}}$, $n\in\mathbb{N}$, satisfies also condition \eqref{con2}, we will use the following approximate functional equation 
\begin{eqnarray}\label{aprfe}
L(s;\chi_j) &=& \sum\limits_{n\leq X}\chi_j(n)n^{-s}+\Delta(s;\chi_j)\sum\limits_{n\leq y}\chi_j^-(n){n^{s-1}}+\nonumber\\
&&+~ O\left(X^{-\sigma}\log(y+2)+\tau^{1/2-\sigma}y^{\sigma-1}\right)
\end{eqnarray}
valid for $s=\sigma+i\tau$ with $0<\sigma<1$, $\tau\geq \tau_0>0$, and $X,y>0$ such that $2\pi Xy=q_j \tau$, where $q_j$ is the modulus of $\chi_j$.
This formula follows from a result due to Vivek Rane \cite{ran}.
Observe also that $L_Q(s;\chi_j)$ can be written as an absolutely convergent Dirichlet series
\begin{equation}\label{trEd}
L_Q(s;\chi_j)=\mathop{\sum\limits_{n=1}^{\infty}}_{p\mid n\Rightarrow p\leq Q}\dfrac{\chi_j(n)}{n^s}
\end{equation}
for any $\sigma>0$ and $Q>0$, as follows from \eqref{trE} by expanding each factor of the truncated Euler product into a geometric series. 

Now let $K$ be a compact subset of the strip $\sigma_1\leq\sigma<1$ for some $\sigma_1\in(1/2,1)$. Then equations \eqref{aprfe} and \eqref{trEd} imply, for every $s=\sigma+it\in K$, any $Q>0$, any sufficiently large $N\in\mathbb{N}$ and for $X=q_jN/(4\pi)$, that
\begin{eqnarray}\label{long}
\lefteqn{\sum\limits_{N\leq n\leq2N}\left|L\left(s+i\gamma_{a}^{{(n)}};\chi_j\right)-L_Q\left(s+i\gamma_{a}^{{(n)}};\chi_j\right)\right|^2}\nonumber\\
&\ll& \sum\limits_{N\leq n\leq2N}\left|\mathop{\sum\nolimits_1}_{m\leq X}\dfrac{\chi_j(m)}{m^{s+i\gamma_{a}^{{(n)}}}}\right|^2+\sum\limits_{N\leq n\leq2N}\left|\mathop{\sum\nolimits_2}_{m> X}\dfrac{\chi_j(m)}{m^{s+i\gamma_{a}^{{(n)}}}}\right|^2\nonumber+\\
&&+ \sum\limits_{N\leq n\leq 2N}\left|\Delta(s+i\gamma_{a}^{{(n)}};\chi_j)\right|^2\left|\sum\limits_{m\leq y}\dfrac{\chi_j^-(m)}{{m^{1-s-i\gamma_{a}^{{(n)}}}}}\right|^2+\nonumber\\
&&+ \sum\limits_{N\leq n\leq 2N}O\left(X^{-2\sigma}\left(\log(y+2)\right)^2+\left(t+\gamma_a^{{(n)}}\right)^{1-2\sigma}y^{2(\sigma-1)}\right)
\end{eqnarray}
where $\sum\nolimits_1$ denotes the sum over integers $m$ for which there is no prime $p\leq Q$ with $p\mid m$, and $\sum\nolimits_2$ denotes the sum over integers $m$ which are divisible only by primes $p\leq Q$. 
We denote the terms on the right-hand side of \eqref{long} by $S_1$, $S_2$, $S_3$, $S_4$, respectively, and we prove that each of them is at most $O\left(NQ^{1-2\sigma_1}\right)$ as $N$ tends to infinity. 
In what follows, asymptotics and limits are not taken with respect to parameter $t$ since $t$ represents the imaginary part of complex numbers $s$ from a compact set $K$. 
The implicit constants depend of course also on $K$ and the finitely many moduli $q_1,\dots,q_J$, but they are negligible in our proof. 
Recall that $\gamma_a^{{(n)}} \sim n/\log n$. Then, we have that
\begin{eqnarray*}
S_4&\ll&\sum\limits_{N\leq n\leq 2N}\left[N^{-2\sigma}\left(\log\left(\dfrac{\gamma_{a}^{{(n)}}}{N}+2\right)\right)^2+\left(\gamma_{a}^{{(n)}}\right)^{1-2\sigma}\left(\dfrac{\gamma_{a}^{{(n)}}}{N}\right)^{2(\sigma-1)}\right]\\
&\ll& N\left[N^{-2\sigma_1}+{N}^{1-2\sigma_1}\log N\right]\\
&=&o(N),
\end{eqnarray*}
while \eqref{pow} implies that
\begin{eqnarray*}
S_3&\ll&\sum\limits_{N\leq n\leq2N}\left(\gamma_{a}^{{(n)}}\right)^{1-2\sigma}\left(\sum\limits_{m\leq y}\dfrac{1}{m^{1-\sigma}}\right)^2\\
&\ll& N\left(\gamma_{a}^{{(N)}}\right)^{1-2\sigma_1}\left(\dfrac{\gamma_{a}^{{(2N)}}}{N}\right)^{2}\\
&\ll& N\left(\dfrac{N}{\log N}\right)^{1-2\sigma_1}\left(\log N\right)^{-2}\\
&=&o(N).
\end{eqnarray*}
To estimate $S_2$, we first observe that its inner sum is a tail of an absolutely convergent series. Therefore,
\begin{eqnarray*}
S_2\ll\sum\limits_{N\leq n\leq2N}\left(\mathop{\sum\nolimits_2}_{m>X}\dfrac{1}{m^{\sigma}}\right)^2\ll N\left(\mathop{\sum\nolimits_2}_{m>N/\log N}\dfrac{1}{m^{\sigma}}\right)^2=o(N).
\end{eqnarray*}
Lastly,
\begin{eqnarray}\label{S_1}
S_1&=&\mathop{\sum\nolimits_1}_{m_1,m_2\leq X}\dfrac{1}{(m_1m_2)^{\sigma}}\left(\dfrac{m_2}{m_1}\right)^{it}\sum\limits_{N\leq n\leq2N}\left(\dfrac{m_2}{m_1}\right)^{i\gamma_{a}^{{(n)}}}\nonumber\\
&\ll&N\mathop{\sum\nolimits_1}_{m\leq X}\dfrac{1}{m^{2\sigma_1}}+\mathop{\sum\nolimits_1}_{1\leq m_1<m_2\leq X}\dfrac{1}{(m_1m_2)^{\sigma_1}}\left|\sum\limits_{N\leq n\leq2N}\left(\dfrac{m_2}{m_1}\right)^{i\gamma_{a}^{{(n)}}}\right|.
\end{eqnarray}
The first term on the right hand-side of \eqref{S_1} is $O\left(NQ^{1-2\sigma_1}\right)$. 
For the second term observe that for any $1\leq m_1<m_2\leq X$, we have $1<m_2/m_1\leq q_jN/(4\pi)$. 
Therefore, Theorem \ref{u.d.} implies that this term is bounded from above by
\begin{eqnarray*}
\mathop{\sum\nolimits_1}_{1\leq m_1<m_2\leq X}\dfrac{\log N}{(m_1m_2)^{\sigma_1}}\left(\left(\log\dfrac{m_2}{m_1}\right)^{-1}+\log\dfrac{m_2}{m_1}\right)&\ll& X^{2-2\sigma_1}\log( X)\log N\\
&\ll&N^{2-2\sigma_1}(\log N)^2\\
&=&o(N).
\end{eqnarray*}

Collecting the above estimates we finally arrive at
$$
\limsup\limits_{N\to\infty}\dfrac{1}{N+1}\sum\limits_{N\leq n\leq2N}\left|L\left(s+i\gamma_{a}^{{(n)}};\chi_j\right)-L_Q\left(s+i\gamma_{a}^{{(n)}};\chi_j\right)\right|^2\ll Q^{1-2\sigma_1}
$$
uniformly in $K$ and arbitrary $Q>0$. Taking $Q$ to infinity shows that the sequence $\gamma_{a}^{{(n)}}$, $n\in\mathbb{N}$, also satisfies \eqref{con2} and thus \eqref{joint} holds.

We prove \eqref{periodic} from \eqref{joint}. We employ techniques introduced by Bagchi \cite{bagchi}, Gonek \cite{gonek0} as well as J\" urgen Sander and the second author \cite{sanste}. Let $\psi\neq\mathbf{0}$ be an $r$-periodic arithmetical function. If $r=1$, then 
$$
L(s;\psi)=\psi(1)\zeta(s)
$$
and \eqref{periodic} holds, since $\zeta(s)$ is $L(s;\chi_0)$, where $\chi_0$ is the Dirichlet character mod $1$ and we apply \eqref{joint} only to this character. 

The case $r=2$ is rather special; here some cases need a restriction on the range of approximation. This observation is due to Jerzy Kaczorowski \cite{jerzy} and has recently been discussed by Philipp Muth and the second author \cite{phil}. Since $\psi(1)\neq\psi(2)$, for any $\sigma>1$, we have
$$
L(s;\psi)=\sum\limits_{n=1}^\infty\dfrac{\psi(1)}{(2n-1)^s}+\sum\limits_{n=1}^\infty\dfrac{\psi(2)}{(2n)^s}=\left(\psi(1)+\dfrac{\psi(2)-\psi(1)}{2^s}\right)\zeta(s)=:P(s)\zeta(s).
$$
The latter holds for all $s\in\mathbb{C}$ by analytic continuation. Observe here that $P(s)$ is analytic and bounded in any half-plane $\sigma\geq\sigma_0$ and 
\begin{equation}\label{Pol}
P(s+i\tau)-P(s)\ll\left\|\tau\dfrac{\log2}{2\pi}\right\|
\end{equation}
uniformly in $\sigma\geq\sigma_0$ and $\tau\in\mathbb{R}$. 
Additionally, $P(s)$ is zero-free in the open set
$$D_0:=\left\{s\in\mathbb{C}:\dfrac{1}{2}<\sigma<1\right\}\setminus\left\{\log\left(1-\dfrac{\psi(2)}{\psi(1)}\right)+2k\pi i:k\in\mathbb{Z}\right\}.$$
Therefore, if we assume that $K\subseteq D_0$ and set
$$
g(s):=\frac{h(s)}{P(s)},
$$
then, from \eqref{joint} we have
\begin{equation*}
\sharp\left\{1\leq n\leq N:\begin{array}{c}\max\limits_{s\in K}\left|\zeta\left(s+i\gamma_a^{{(n)}}\right)-g(s)\right|<\eta,\\ \&\ \left\|\gamma_a^{{(n)}}\dfrac{\log 2}{2\pi}\right\|<\eta
\end{array}\right\}>c(\eta)N
\end{equation*}
for any $\eta>0$ and any sufficiently large $N\gg_{\eta}1$, where $c(\eta)>0$ is constant. For those $n$ from the set described on the left-hand side above, in combination with \eqref{Pol}, it also follows that
\begin{eqnarray*}
\max\limits_{s\in K}\left|L\left(s+i\gamma_a^{{(n)}};\psi\right)-h(s)\right|&<&\max\limits_{s\in K}\hspace{-0,5pt}\left|P\hspace{-0,5pt}\left(s\hspace{-0,5pt}+\hspace{-0,5pt}i\gamma_a^{{(n)}}\right)\hspace{-0,5pt}\right|\hspace{-0,5pt}\max\limits_{s\in K}\hspace{-0,5pt}\left|\zeta\hspace{-0,5pt}\left(s\hspace{-0,5pt}+\hspace{-0,5pt}i\gamma_a^{{(n)}}\right)\hspace{-0,5pt}-g(s)\right|\\
&&+~ \max\limits_{s\in K}|g(s)|\max\limits_{s\in K}\left|P\left(s+i\gamma_a^{{(n)}}\right)-P(s)\right|\\
&\ll&\eta.
\end{eqnarray*}
Taking $0<\eta\ll\varepsilon$ sufficiently small, we obtain \eqref{periodic} with the restriction we imposed on $K$.

If $r\geq3$, then $\phi(r)\geq2$, where $\phi$ is the Euler totient function. 
Assuming that $\psi$ is a multiple of a Dirichlet character mod $r$, we can work as in the case of $r=1$. If $\psi$ is not such a multiple, then for every $s\in\mathbb{C}$
\begin{eqnarray}\label{licom1}
L(s;\psi)&=&\dfrac{1}{r^s}\sum\limits_{n=1}^r\psi(n)\zeta\left(s;\dfrac{n}{r}\right)\nonumber\\
&=&\dfrac{1}{r^s}\sum\limits_{n=1}^r\psi(n)\dfrac{r^s}{\phi(r)}\sum\limits_{i=1}^{\phi(r)}\overline{\chi_i}(n)L(s;\chi_i)\nonumber\\
&=&\sum\limits_{i=1}^{\phi(r)}\left(\dfrac{1}{\phi(r)}\sum\limits_{n=1}^r\psi(n)\overline{\chi_i}(n) \right)L(s;\chi_i),
\end{eqnarray}
where $\chi_i$, $i=1,2,\dots,\phi(r)$, are the Dirichlet characters mod $r$.
The expression of $L(s;\psi)$ as a linear combination of Hurwitz zeta-functions with rational parameters follows first for $\sigma>1$, where there are absolutely convergent Dirichlet series representations of them, and then by analytic continuation to the whole complex plane.
The expression of a Hurwitz zeta-function with rational parameter as a linear combination of Dirichlet $L$-functions follows from the orthogonality relation of the characters.
Now if we set
\begin{eqnarray*}
c_i:=\dfrac{1}{\phi(r)}\sum\limits_{n=1}^r\psi(n)\overline{\chi_i}(n),&i=1,2,\dots,\phi(r),
\end{eqnarray*}
then at least two of them, say $c_1$ and $c_2$, are non-zero by assumption.
If we define $M_h:=1+\max_{s\in K}|h(s)|$ and the functions 
\begin{eqnarray}\label{licom2}
g_1(s):=\dfrac{h(s)+M_h}{c_1},&g_2(s):=-\dfrac{M_h}{c_2},
\end{eqnarray}
\begin{eqnarray}\label{licom3}
g_i(s):=\eta,&i=3,\dots,\phi(r)
\end{eqnarray}
for a given $\eta>0$, then $g_i$, $i=1,2,\dots,\phi(r)$, are non-zero continuous functions on $K$ which are analytic in its interior and
\begin{equation}\label{h(s)}
h(s)=\sum\limits_{i=1}^{\phi(r)}c_ig_i(s)-\sum\limits_{i=3}^{\phi(r)}c_ig_i(s)=\sum\limits_{i=1}^{\phi(r)}c_ig_i(s)-\eta(\phi(r)-3)\sum\limits_{i=3}^{\phi(r)}c_i.
\end{equation}
Since $\chi_i$, $i=1,2,\dots,\phi(r)$, are non-equivalent Dirichlet characters, by \eqref{joint} we obtain
\begin{equation*}
\sharp\left\{1\leq n\leq N:\max\limits_{1\leq i\leq\phi(r)}\max\limits_{s\in K}\left|L\left(s+i\gamma_a^{{(n)}};\chi_i\right)-g_i(s)\right|<\eta\right\}>c(\eta)N
\end{equation*}
for any sufficiently large $N\gg_{\eta}1$.
In this case, for those $n$ from the set described on the left-hand side above, in combination with \eqref{licom1}--\eqref{licom3}, it also follows that
\begin{eqnarray*}
\max\limits_{s\in K}\left|L\left(s+i\gamma_a^{{(n)}};\psi\right)-h(s)\right| &\leq& \sum\limits_{i=1}^{\phi(r)}|c_i|\max\limits_{s\in K}\left|L\left(s+i\gamma_a^{{(n)}};\chi_i\right)-g_i(s)\right| \\
&&+~ \eta(\phi(r)-3)\sum\limits_{i=3}^{\phi(r)}|c_i|\\
&\ll& \eta.
\end{eqnarray*}
Taking $0<\eta\ll\varepsilon$ sufficiently small, we obtain \eqref{periodic}. Observe that in this case, $h(s)$ is allowed to have zeros.

\bigskip


\bigskip

\noindent
Athanasios Sourmelidis\\
Institute of Analysis and Number Theory, TU Graz\\
Steyrergasse 30, 8010, Graz, Austria\\
sourmelidis@math.tugraz.at
\medskip

\noindent
J\"orn Steuding\\
Department of Mathematics, W\"urzburg University\\
Emil Fischer-Str. 40, 97\,074 W\"urzburg, Germany\\
steuding@mathematik.uni-wuerzburg.de
\medskip

\noindent
Ade Irma Suriajaya\\
Faculty of Mathematics, Kyushu University \\
744 Motooka, Nishi-ku, Fukuoka 819-0395, Japan\\
adeirmasuriajaya@math.kyushu-u.ac.jp

\end{document}